\documentclass[11pt]{article}
\usepackage{mathrsfs}
\usepackage{amssymb}
\usepackage{amsfonts}
\usepackage{color,xcolor}
\usepackage{graphicx}
\usepackage{geometry}
\usepackage{manfnt}
\usepackage[pagewise]{lineno}
\usepackage[colorlinks,citecolor=blue,urlcolor=blue]{hyperref}
\newtheorem{theorem}{Theorem}[section]
\newtheorem{corollary}{Corollary}[section]
\newtheorem{lemma}{Lemma}[section]

\newtheorem{definition}{Definition}[section]
\newtheorem{remark}{Remark}[section]
\newtheorem{example}{Example}[section]

\def\[{{\Big[}}\def\]{{\Big]}}\def\({{\Big(}}\def\){{\Big)}}

\def\cC{{\mathcal C}}

\def\mE{{\mathbb E}}

\def\mN{{\mathbb N}}\def\mP{{\mathbb P}}
\def\mR{{\mathbb R}}

\def\={&\!\!=\!\!&}

\def\geq{\geqslant}\def\leq{\leqslant}

\def\div{\mathord{{\rm div}}}

\begin{document}
\title{ \bf Stochastic transport equation with bounded and Dini continuous drift}

\author{Jinlong Wei$^a$, Guangying Lv$^b$ and Wei Wang$^c$
\\ {\small \it $^a$School of Statistics and Mathematics, Zhongnan University of}\\
{\small \it Economics and Law, Wuhan 430073, China}
\\ {\small \tt  weijinlong.hust@gmail.com}
\\ {\small \it $^b$College of Mathematics and Statistics,
Nanjing University of Information}
\\{\small \it  Science and Technology, Nanjing 210044, China}\\ {\small \tt  gylvmaths@126.com}
\\ {\small \it $^c$Department of Mathematics, Nanjing University, Nanjing 210093, China}
\\ {\small \tt wangweinju@nju.edu.cn} }
\date{\today}

 \maketitle
\noindent{\hrulefill}

\begin{abstract}
The results established by Flandoli, Gubinelli
and Priola ({\it Invent. Math.} {\bf 180} (2010) 1--53) for stochastic transport equation with bounded and H\"{o}lder continuous drift are generalized  to bounded and Dini continuous drift. The  uniqueness of $L^\infty$-solutions  is established by the It\^o--Tanaka trick partially solving the uniqueness problem, which is still open,   for  stochastic transport equation with  only bounded measurable drift.    Moreover    the existence and uniqueness of stochastic diffeomorphisms flows for a stochastic differential
equation with bounded and Dini continuous drift is obtained.
\end{abstract}


 \vskip0mm\noindent{\hrulefill}

\section{Introduction}\label{sec1}\setcounter{equation}{0}
We are concerned with the following stochastic transport equation
\begin{eqnarray}\label{1.1}
\left\{
  \begin{array}{ll}
\partial_tu(t,x)+b(t,x)\cdot\nabla u(t,x)
+\sum_{i=1}^d\partial_{x_i}u(t,x)\circ\dot{B}_i(t)=0, \quad (t,x)\in(0,T)\times {\mathbb R}^d, \\
u(t,x)|_{t=0}=u_0(x), \quad  x\in{\mathbb R}^d,
  \end{array}
\right.
\end{eqnarray}
where $\{B(t)\}_t=\{(B_1(t), B_2(t), _{\cdots}, B_d(t))\}_t$ is a
$d$-dimensional standard Brownian motion defined on a stochastic
basis ($\Omega, {\mathcal F},{\mathbb P},({\mathcal F}_{t})_{t\geq 0}$).  The stochastic
integration with  notation $\circ$ is interpreted in Stratonovich
sense. Given $T>0$, the drift coefficient  $b: [0,T]\times{\mathbb R}^d\rightarrow{\mathbb R}^d$ and the initial
data~$u_0:
{\mathbb R}^d\rightarrow{\mathbb R}$ are measurable functions  in $L^1([0,T];L^1_{loc}({\mathbb R}^d;{\mathbb R}^d))$
and $L^1([0,T];L^1_{loc}({\mathbb R}^d))$ respectively.  We are interested in the existence and uniqueness of weak
$L^\infty$-solutions for the stochastic equation~(\ref{1.1}).

There are many  results on the weak $L^\infty$-solutions for the deterministic transport equation.    The first remarkable result   of the uniqueness  solution in $L^{\infty}([0, T]\times\mathbb{R}^{d})$ was obtained by  DiPerna and Lions~\cite{DL}
under the assumption $b\in L^1([0,T];W^{1,1}_{loc}({\mathbb R}^d;{\mathbb R}^d))$ with suitable global conditions
including $L^\infty$-bounds on spatial divergence.  Then Ambrosio \cite{Amb} weakened the condition $W^{1,1}_{loc}$ to~$BV_{loc}$ and the uniqueness of weak $L^\infty$-solutions  is obtained under the assumption that negative part of $\div b$ is in space  $L^1([0,T];L^\infty(\mR^d))$.  For general $b$ just with $BV_{loc}$ or H\"{o}lder regularity,  the uniqueness of weak $L^\infty$-solutions for the deterministic
equation  fails and counterexamples    have been constructed  in many works~\cite{Bre,CLR1,CD,Dep,DL}.
Obviously, more restrictions need to be imposed on $b$ to overcome the obstacle of nonuniqueness of solution in deterministic case.

However the appearance of the noise makes the solution unique under very general assumptions on the drift coefficient for ordinary differential equations \cite{KR,Ver}, so a natural idea is to investigate the effects of noise  in transport equation. The first milestone result, founded by Flandoli, Gubinelli and Priola~\cite{FGP1}, showed the uniqueness of weak $L^\infty$-solutions just with assuming~$b\in L^\infty([0,T];\cC_b^\alpha(\mR^d;\mR^d))$ and $\div b\in L^p([0,T]\times\mR^d)$ for some $\alpha>0$, $p>2$.
This is the first concrete example of a partial differential equation related to fluid dynamics that  becomes
well-posed with a suitable noise. A key step  for  this result is to  perform differential computations on regularization of $L^\infty$-solutions by a  commutator lemma.  Unfortunately  the strategy fails  if $b$ is not Sobolev differentiable.  However,   by observing the fact that   stochastic differential equation~(SDE)
\begin{eqnarray}\label{1.2}
\left\{\begin{array}{ll}\medskip
dX(t)=b(t,X(t))dt+dB(t), \quad  0<t\leq T,\\
X(t)|_{t=0}=x\in{\mathbb R}^d,  \end{array}\right.
\end{eqnarray}
defines  a $\cC^1$ stochastic diffeomorphisms flow and,  along the stochastic characteristic $X(t)$, the integral
$\int_0^t\div b(s,X(s,x))ds$ has a regularization,   Flandoli, Gubinelli and Priola developed
the commutator lemma to prove the uniqueness of  solutions.
On the other hand,  for bounded measurable $b$,  Mohammed,  Nilssen and Proske~\cite{MNP} also proved the existence, uniqueness and Sobolev differentiable stochastic
flows for~(\ref{1.2}) by employing ideas from  the Malliavin
calculus coupled with new probabilistic estimates on the spatial weak derivatives of solutions of~(\ref{1.2}).
Then, as an application, they obtained the existence and uniqueness of Sobolev differentiable
weak solutions for (\ref{1.1}) with every $\cC^1(\mR^d)$ initial data.  Notice that in one result  the stochastic flow $\{X(t, x)\}$ is  differentiable in $x$ in the classical sense~\cite{FGP1} and in  the other result the stochastic flow $\{X(t,x)\}$ is only Sobolev differentiable~\cite{MNP}. So the method developed by Mohammed,  Nilssen and Proske~\cite{MNP} can not be adapted to establish the uniqueness of weak $L^\infty$-solutions to (\ref{1.1}).

Recent result~\cite{AF}, by using a different philosophy, proved  the uniqueness of weak $L^\infty$-solutions for~(\ref{1.1})  just with  assuming the $BV$ regularity for $b$ but without the $L^\infty$-bounds on spatial divergence. There are also several other related works~\cite{BFG,FL,FF2,NO,WDGL,Zhang1}\,.
  However, \textbf{for bounded measurable $b$ and $\div b\in L^p([0,T]\times \mR^d)$ (for some $p\in [1,\infty)$), the uniqueness of weak
$L^\infty$-solutions for (\ref{1.1}) is still unknown.}\\

This paper intend  to give a partial answer for the above problem and  novelties  of the work  are
 \begin{itemize}
\item   {\it The uniqueness of weak $L^\infty$-solutions for the Cauchy problem (\ref{1.1}) with bounded and
Dini-continuous drift is established due to the existence of noise, while  the  corresponding deterministic equation has  multiple solutions.}
\item  {\it The
existence and uniqueness of stochastic diffeomorphisms flow for singular SDE (\ref{3.1}) is established without H\"{o}lder continuity
or  Sobolev differentiability hypotheses on $b$.}

\item  {\it The maximum regularity for parabolic equations of second order with  H\"{o}lder-Dini or strong H\"{o}lder or weak H\"{o}lder coefficients is established.}
\end{itemize}

We follow the strategy of Flandoli, Gubinelli and Priola's~\cite{FGP1} to  establish the existence of  a
stochastic~$\cC^1$ diffeomorphisms flow for (\ref{3.1}) by  the It\^o--Tanaka trick, then  derive a commutator
estimates to get the uniqueness for weak $L^\infty$-solution of (\ref{1.1}).  The main idea of It\^o--Tanaka trick is to use a  parabolic partial differential equation~(PDE) to transform the original SDE~(\ref{3.1}) with irregular
drift and regular diffusion to  a new SDE~(\ref{3.15}) with regular drift and diffusion.  Then by the equivalence between (\ref{3.1}) and (\ref{3.15}) we show the existence of the  stochastic $\cC^1$ diffeomorphisms flow for  SDE (\ref{3.1})\,.
There are also some recent works on  the stochastic flows and SDEs~\cite{Att,FF1,FF3,FGP2,TWT,WZ,Zhang2}.

\medskip
In the following parts, we first derive the $W^{2,\infty}$ estimates for a class of
second order parabolic PDEs with bounded and Dini continuous coefficients in section 2; then by using the  $W^{2,\infty}$ estimates,  the existence and uniqueness of stochastic flow of diffeomorphisms  for
SDE (\ref{3.1}) is shown in section~3 by the It\^o--Tanaka trick; last section is concerned with the existence and uniqueness of weak $L^\infty$-solutions to stochastic transport equation (\ref{1.1}).

\vskip2mm\noindent
\textbf{Notations} The letter $C$  denotes a positive constant, whose values may change in different places. For a
parameter or a function $\kappa$, $C(\kappa)$ means the constant is only dependent on $\kappa$, and we also write it
as $C$ if there is no confusion. ${\mathbb N}$ is the set of natural numbers. For every
$R>0$, $B_R:=\{x\in{\mathbb R}^d:|x|<R\}$. Almost surely is  abbreviated to $a.s.$.
Let $\Theta$ be a ${\mathbb R}^{d\times d}$-valued function $\Theta=(\Theta_{i,j}(x))_{d\times d}$ with norm
$\|\Theta(x)\|=\max_{1\leq i,j\leq d}|\Theta_{i,j}(x)|$. For $\xi\in \mR^d$,   $|\xi|=(\sum_{i=1}^d\xi_i^2)^{1/2}$. ${\mathbb R}_+$ is the set of nonnegative real numbers and $\bar{{\mathbb R}}_+={\mathbb R}_+\cup\{+\infty\}$.

\section{Parabolic PDEs with  bounded and Dini coefficients}\label{sec2}\setcounter{equation}{0}
Let $T>0$. Consider the following Cauchy problem
\begin{eqnarray}\label{2.1}
\left\{\begin{array}{ll}\medskip
\partial_{t}u(t,x)=\frac{1}{2}\Delta u(t,x)+g(t,x)\cdot \nabla u(t,x)+f(t,x), \ (t,x)\in (0,T)\times {\mathbb R}^d, \\
u(0,x)=0, \  x\in{\mathbb R}^d.  \end{array}\right.
\end{eqnarray}
 The function $u(t,x)$ is called a strong solution of (\ref{2.1}) if $ u\in L^\infty([0,T];W^{2,\infty}({\mathbb R}^d))\cap W^{1,\infty}([0,T];L^\infty({\mathbb R}^d))$ such that for almost all $(t,x)\in [0,T]\times\mR^d$, (\ref{2.1}) holds. We have the following  equivalent form for the strong solution.

\begin{lemma} \label{lem2.1}
Let $f\in L^\infty([0,T];\cC_b(\mR^d))$,  $g\in L^\infty([0,T];\cC_b(\mR^d;\mR^d))$ and  $u\in L^\infty([0,T];W^{2,\infty}({\mathbb R}^d))\cap W^{1,\infty}([0,T];L^\infty({\mathbb R}^d))$, then $u$ is a strong solution for (\ref{2.1}) if and only if
\begin{eqnarray}\label{2.2}
u(t,x)&=&\int\limits_0^tK(t-s,\cdot)\ast (g(s,\cdot)\cdot \nabla
u(s,\cdot))(x)ds\nonumber \\ && +\int\limits_0^tK(t-s,\cdot)\ast f(s,\cdot)(x)ds,  \quad  for\; all  \ (t,x)\in [0,T]\times\mR^d,
\end{eqnarray}
where $K(t,x)=(2\pi t)^{-\frac{d}{2}}e^{-\frac{|x|^2}{2t}}, \ t>0, \ x\in\mathbb{R}^d$.
\end{lemma}

\vskip2mm\noindent\textbf{Proof.}
By the properties of  the heat kernel $K$, it is direct to verify that  if $u$ satisfies (\ref{2.2}),  for almost all $(t,x)\in [0,T]\times\mR^d$, (\ref{2.1}) holds. On the other hand, if $u$ satisfies (\ref{2.1}) then for every $\psi\in \mathcal{C}_0^\infty({\mathbb R}^d)$ and every $t\in [0,T]$
\begin{eqnarray*}
\int\limits_0^t\int\limits_{\mathbb{R}^d}\partial_su(s,x)\varphi(s,x)dsdx&=&
\frac{1}{2}\int\limits_{\mathbb{R}^d}\int\limits_0^t\Delta u(s,x)\varphi(s,x)dsdx+ \int\limits_{\mathbb{R}^d}\int\limits_0^t
g(s,x)\cdot \nabla u(s,x) \varphi(s,x) dsdx \nonumber\\&&+\int\limits_{\mathbb{R}^d}\int\limits_0^t f(s,x)\varphi(s,x)dsdx,
\end{eqnarray*}
with  $\varphi(s,x)=K(t-s,\cdot)\ast \psi(\cdot)(x)$.

Now integrating  by parts yields
\begin{eqnarray*}
\int\limits_{\mathbb{R}^d}u(t,x)\psi(x)dx&=&
\int\limits_{\mathbb{R}^d}\int\limits_0^t u(s,x)[\partial_s\varphi(s,x)+\frac{1}{2}\Delta\varphi(s,x)]dsdx \nonumber\\&&+ \int\limits_{\mathbb{R}^d}\int\limits_0^t
g(s,x)\cdot \nabla u(s,x) \varphi(s,x) dsdx+\int\limits_{\mathbb{R}^d}\int\limits_0^t f(s,x)\varphi(s,x)dsdx
\nonumber\\ &=& \int\limits_{\mathbb{R}^d}\int\limits_0^tK(t-s,\cdot)\ast (g(s,\cdot)\cdot \nabla u(s,\cdot))(x)ds\psi(x)dx\nonumber\\&&+\int\limits_{\mathbb{R}^d}\int\limits_0^tK(t-s,\cdot)\ast f(s,\cdot)(x)ds\psi(x)dx, \quad  for\; all  \ t\in [0,T],
\end{eqnarray*}
then by the arbitrariness of  $\psi$   and  continuity of $u$  in $x$, (\ref{2.2}) holds. $\Box$

\begin{definition} \label{def2.1} An increasing continuous  function $\phi: \mR_+\rightarrow \mR_+$  is called  a Dini function if
\begin{eqnarray}\label{2.3}
\int\limits_{0+}\frac{\phi(r)}{r}dr<+\infty.
\end{eqnarray}
A measurable function $h:\mR^{d}\rightarrow \mR$  is said to be Dini continuous if there is a Dini function $\phi$ such that
\begin{eqnarray}\label{2.4}
|h(x)-h(y)|\leq \phi(|x-y|).
\end{eqnarray}
\end{definition}

We now state the  main result of this section.
\begin{theorem} \label{the2.1}  Let $f\in L^\infty([0,T];\cC_b(\mR^d))$ and $g\in L^\infty([0,T];\cC_b(\mR^d;\mR^d))$. Suppose that $r_0\in (0,1)$ and  there is a Dini function $\phi$ such that for every $x\in\mR^d$
\begin{eqnarray}\label{2.5}
|f(t,x)-f(t,y)|+|g(t,x)-g(t,y)|\leq \phi(|x-y|), \quad  for\; all  \ y\in B_{r_0}(x), \ t\in [0,T].
\end{eqnarray}

(i) The Cauchy problem (\ref{2.1}) has a unique strong solution $u$, and  there is a constant $C(d,T)$ such that
\begin{eqnarray}\label{2.6}
\|u\|_{L^\infty([0,T];\cC^2_b({\mathbb R}^d))} \leq C(d,T)(1+\|f\|_{L^\infty([0,T];\cC_b(\mR^d))}+
\|g\|_{L^\infty([0,T];\cC_b(\mR^d;\mR^d))}).
\end{eqnarray}
Moreover, for every $1\leq i,j\leq d$ and every $x,y\in \mR^d$, there is another constant $C(d,T)$ such that
\begin{eqnarray}\label{2.7}
&&|\partial^2_{x_i,x_j}u(t,x)-\partial^2_{y_i,y_j}u(t,y)| \nonumber \\ &\leq &C(d,T)\left[\int\limits_{r\leq |x-y|}\frac{\phi(r)}{r}dr+ \phi(|x-y|)+|x-y|\int\limits_{|x-y|<r\leq r_0}\frac{\phi(r)}{r^2}dr\right]1_{|x-y|<r_0}\nonumber \\ &&+C(d,T)|x-y|, \quad    for \;  all  \ t\in [0,T].
\end{eqnarray}

(ii)  Let $\varrho_n$, $n\in \mN$,  be a regularizing kernel, that is
\begin{eqnarray}\label{2.8}
\varrho_n(x) =n^d \varrho(nx) \ \ \mbox{with} \ \ 0\leq \varrho \in \mathcal{C}^\infty_0({\mathbb R}^d) , \ \ \mbox{supp}(\varrho)\subset B_1, \ \ \int\limits_{{\mathbb R}^d}\varrho(x)dx=1.
\end{eqnarray}
  Let $u^n$ be the unique strong solution of (\ref{2.1}) with $f$ and $g$ are replaced by  $f^n(t,x)=f\ast \varrho_n(t,x)$ and $g^n(t,x)=g\ast \varrho_n(t,x)$ respectively. Then $u^n\in L^\infty([0,T];\cC^2_b({\mathbb R}^d))\cap W^{1,\infty}([0,T];\cC_b({\mathbb R}^d))$ and satisfies~(\ref{2.6})--(\ref{2.7}) uniformly in $n$. Furthermore,
\begin{eqnarray}\label{2.9}
\lim_{n\rightarrow\infty}\| u^n-u\|_{L^\infty([0,T];\cC_b^2({\mathbb R}^d))}=0.
\end{eqnarray}
\end{theorem}
\noindent\textbf{Proof.} (i) We first prove the result  for the case  $g=0$.  By Lemma \ref{lem2.1}  we just need   to show
\begin{eqnarray}\label{2.10}
u(t,x)=\int\limits_0^{t}K(t-s,\cdot)\ast f(s,\cdot)(x)ds
\end{eqnarray}
is in $L^\infty([0,T];\cC^{2}_b({\mathbb R}^d))$ and (\ref{2.7}) holds. In fact   $u\in L^\infty([0,T];W^{1,\infty}({\mathbb R}^d))$ is classical~\cite[Ch.4]{LSU} by the explicit representation (\ref{2.10})\,.

 Next  we    show $\partial^2_{x_i,x_j}u\in L^\infty([0,T];\cC_b({\mathbb R}^d))$ for every $1\leq i,j\leq d$ and (\ref{2.7}) holds. For this we first show that  $\partial^2_{x_i,x_j}u\in L^\infty([0,T];L^\infty({\mathbb R}^d))$.
Let $\theta\in (0,1/2)$. For $x\in \mR^d$ and $t\in (0,T]$, we have
\begin{eqnarray}\label{2.11}
\Big|\partial^2_{x_i,x_j}u(t,x)\Big|&=&\Bigg|\int\limits_0^t\int\limits_{{\mathbb R}^d}\partial^2_{x_i,x_j}K(t-s,x-y)f(s,y)dyds\Bigg|\nonumber\\
&=& \Bigg|\int\limits_0^t\int\limits_{{\mathbb R}^d}\partial^2_{x_i,x_j}K(t-s,x-y)[f(s,y)-f(s,x)]dyds \Bigg|
\nonumber\\
&=& \Bigg|\int\limits_0^t\int\limits_{|x-y|> (t-s)^\theta }\partial^2_{x_i,x_j}K(t-s,x-y)[f(s,y)-f(s,x)]dyds\nonumber\\&&+\int\limits_0^t\int\limits_{|x-y|\leq (t-s)^\theta }\partial^2_{x_i,x_j}K(t-s,x-y)[f(s,y)-f(s,x)]dyds\Bigg|
\nonumber\\
&\leq & C\|f\|_{L^\infty([0,T]\times\mR^d)}\int\limits_0^t\int\limits_{|x-y|> (t-s)^\theta }\frac{1}{t-s}K(t-s,x-y)dyds\nonumber\\&&+\int\limits_0^t\int\limits_{|x-y|<(t-s)^\theta \wedge r_0 }\frac{1}{t-s}K(t-s,x-y)\phi(|x-y|)dyds \nonumber\\
&&+C\int\limits_0^t\int\limits_{(t-s)^\theta \wedge r_0\leq |x-y|\leq(t-s)^\theta }\frac{1}{t-s}K(t-s,x-y)|x-y|dyds
\nonumber\\&\leq & C\[\int\limits_0^T\int\limits_{|y|> s^\theta }\frac{1}{s}K(s,y)dyds+\int\limits_0^T\int\limits_{|y|<s^\theta \wedge r_0 }\frac{1}{s}K(s,y)\phi(|y|)dyds
\nonumber\\
&&+C\int\limits_{r_0^{\frac{1}{\theta}}}^T\int\limits_{r_0\leq |y|\leq s^\theta }\frac{1}{s}K(s,y)|y|dyds\]
=:I_1+I_2+I_3,
\end{eqnarray}
where the boundedness of $f$ is applied  in the seventh line of (\ref{2.11}), and in the last line   $I_3=0$ for~$T\leq r_0^{\frac{1}{\theta}}$.
We first estimate $I_1$.
\begin{eqnarray}\label{2.12}
I_1&=& C\int\limits_0^Ts^{-\frac{d+2}{2}}ds\int\limits_{r> s^\theta}e^{-\frac{r^2}{2s}}r^{d-1}dr
\nonumber\\&=& C\int\limits_0^{T}s^{-1}ds\int\limits_{s^{\theta-\frac12}}^{\infty} e^{-\frac{r^2}{2}}r^{d-1}dr
\nonumber\\&=&C\int\limits_0^{\frac12}s^{-1}ds\int\limits_{s^{\theta-\frac12}}^{\infty} e^{-\frac{r^2}{2}}r^{d-1}dr+C\int\limits_{\frac12}^{T\vee {\frac12}}s^{-1}ds\int\limits_{s^{\theta-\frac12}}^{\infty} e^{-\frac{r^2}{2}}r^{d-1}dr
\nonumber\\&\leq &C\int\limits_0^{\frac12}s^{-1}ds\int\limits_{s^{\theta-\frac12}}^{\infty} e^{-\frac{r^2}{2}}r^{d-1}dr+C\log(2T\vee 1)\int\limits_{(T\vee \frac12)^{{\theta-\frac12}}}^{\infty} e^{-\frac{r^2}{2}}r^{d-1}dr
\nonumber\\&\leq &C\Bigg[\log(2T\vee 1)+\int\limits_0^{\frac12}s^{-1}ds\int\limits_{s^{\theta-\frac12}}^{\infty} e^{-\frac{r^2}{2}}r^{d-1}dr\Bigg].
\end{eqnarray}
Without loss of generality we assume that $d$ is even.
Otherwise, since $\theta<1/2$ and $s\in (0,1/2)$,
\begin{eqnarray*}
\int\limits_{s^{\theta-\frac12}}^{\infty} e^{-\frac{r^2}{2}}r^{d-1}dr\leq \int\limits_{s^{\theta-\frac12}}^{\infty} e^{-\frac{r^2}{2}}r^ddr.
\end{eqnarray*}
  Thus, there is a natural number $m\geq 0$ such that $d=2m+2$ and
\begin{eqnarray}\label{2.13}
\int\limits_{s^{\theta-\frac12}}^{\infty} e^{-\frac{r^2}{2}}r^{d-1}dr=2^m \int\limits_{s^{2\theta-1}/2}^{\infty} e^{-r}r^mdr=:2^mJ_m.
\end{eqnarray}
Set $s_0=s^{2\theta-1}/2$ and integrating  by parts yields the  following recurrence formula
\begin{eqnarray*}
J_m=\int\limits_{s^{2\theta-1}/2}^{\infty} e^{-r}r^mdr=s_0^me^{-s_0}+mJ_{m-1}.
\end{eqnarray*}
Then
\begin{eqnarray}\label{2.14}
J_m\leq C(1+s_0^m)e^{-s_0}.
\end{eqnarray}
Now  from (\ref{2.12})--(\ref{2.14})
\begin{eqnarray}\label{2.15}
I_1&\leq& C\Bigg[\log(2T\vee 1)+\int\limits_0^{\frac12}s^{-1}(1+s^{(\theta-\frac12)(d-2)})
e^{-\frac{s^{2\theta -1}}{2}}
ds\Bigg]
\nonumber\\&
=&C\Bigg[\log(2T\vee 1)+\int\limits_2^{\infty} s^{-1}(1+s^{(\frac12-\theta)(d-2)})
e^{-\frac{s^{1-2\theta}}{2}}ds\Bigg]<+\infty.
\end{eqnarray}

For $I_{2}$\,,  since $\phi$ is a nonnegative increasing continuous function, we have
\begin{eqnarray}\label{2.16}
I_2\leq C\int\limits_0^{T}\frac{\phi(s^\theta)}{s}ds\int\limits_{\mR^d} K(s,y)dy =C\int\limits_0^{T^\theta}\frac{\phi(s)}{s}ds<+\infty.
\end{eqnarray}
Last for $I_3$, we have
\begin{eqnarray}\label{2.17}
I_3\leq C\int\limits_{0}^T\int\limits_{\mR^d}s^{-\frac{1}{2}}K(s,y)dyds\leq C\sqrt{T}< +\infty.
\end{eqnarray}
Now by (\ref{2.11}) and  (\ref{2.15})--(\ref{2.17}),    $u\in L^\infty([0,T];W^{2,\infty}({\mathbb R}^d))$ and there is a constant $C>0$\,,  depending  only on $d$ and $T$\,,  such that
\begin{eqnarray*}
\|u\|_{L^\infty([0,T];W^{2,\infty}({\mathbb R}^d))} \leq C(1+\|f\|_{L^\infty([0,T];\cC_b(\mR^d))}).
\end{eqnarray*}
So  to prove $u\in L^\infty([0,T];\cC^{2}_b({\mathbb R}^d))$,  we just need inequality (\ref{2.7}).  Notice that  for $x,y\in\mR^d$ and $|x-y|\geq r_0/2$, since  $u\in L^\infty([0,T];W^{2,\infty}({\mathbb R}^d))$, there is constant $C(d,T)$ such that
\begin{eqnarray*}
|\partial^2_{x_i,x_j}u(t,x)-\partial^2_{y_i,y_j}u(t,y)|\leq \frac{C(d,T)r_0}{2} \leq C(d,T)|x-y|, \quad {\rm for \; all} \ \ t\in [0,T].
\end{eqnarray*}
We next show (\ref{2.7}) for $0<|x-y|<r_0/2$.

 For every $1\leq i,j\leq d$ and $0<t\leq T$, every $x,y\in \mR^d$  with $|x-y|\leq r_0/2$, from  (\ref{2.10})\,,  we have
\begin{eqnarray}\label{2.18}
&&\partial^2_{x_i,x_j}u(t,x)-\partial^2_{y_i,y_j}u(t,y)\nonumber\\
&=&\int\limits_0^t
ds\int\limits_{|x-z|\leq 2|x-y|}\partial^2_{x_i,x_j}K(t-s,x-z)[f(s,z)-f(s,x)]dz
\nonumber\\&&-\int\limits_0^t
ds\int\limits_{|x-z|\leq 2|x-y|}\partial^2_{y_i,y_j}K(t-s,y-z)[f(s,z)-f(s,y)]dz
\nonumber\\&&+\int\limits_0^t
ds\int\limits_{|x-z|> 2|x-y|}\partial^2_{y_i,y_j}K(t-s,y-z)[f(s,y)-f(s,x)]dz
\nonumber\\
&&+\int\limits_0^t
ds\int\limits_{|x-z|> 2|x-y|}[\partial^2_{x_i,x_j}K(t-s,x-z)-\partial^2_{y_i,y_j}K(t-s,y-z)][f(s,z)-f(s,x)]dz
\nonumber\\ &=&:J_1(t)+J_2(t)+J_3(t)+J_4(t).
\end{eqnarray}

By the assumption (\ref{2.5})
\begin{eqnarray}\label{2.19}
|J_1(t)|&\leq&  \int\limits_0^tds\int\limits_{|x-z|\leq 2|x-y|} s^{-\frac{d+2}{2}}e^{-\frac{|x-z|^2}{2s}}\phi(|x-z|) dz
\nonumber\\&=& \int\limits_{|z|\leq 2|x-y|} \phi(|z|)dz  \int\limits_0^t s^{-\frac{d+2}{2}}e^{-\frac{|z|^2}{2s}}ds
\nonumber\\ &\leq&C \int\limits_{|z|\leq 2|x-y|} \frac{\phi(|z|)}{|z|^d} dz  \int\limits_0^\infty r^{\frac{d}{2}-1}e^{-\frac{r}{2}}dr
\nonumber\\&\leq&C\int\limits_{r\leq 2|x-y|} \frac{\phi(r)}{r} dr.
\end{eqnarray}

Analogously,
\begin{eqnarray}\label{2.20}
|J_2(t)|\leq  C\int\limits_{r\leq 2|x-y|} \frac{\phi(r)}{r} dr.
\end{eqnarray}

For $J_3$,  by Gauss--Green's formula
\begin{eqnarray}\label{2.21}
|J_3(t)|&=&\left|\int\limits_0^t
ds\int\limits_{|x-z|=2|x-y|}\partial_{y_j}K(t-s,y-z)n_i[f(s,y)-f(s,x)]dS\right|
\nonumber\\
 &\leq& C\int\limits_0^{T}ds
\int\limits_{|x-z|=2|x-y|}|y-z|s^{-\frac{d+2}{2}}e^{-\frac{|y-z|^2}{2s}}\phi(|x-y|)dS
\nonumber\\&\leq&
C\phi(|x-y|)\int\limits_{|x-z|=2|x-y|}|y-z| dS \int\limits_0^{T}
s^{-\frac{d+2}{2}}e^{-\frac{|x-y|^2}{2s}}ds
\nonumber\\ &\leq& C\phi(|x-y|)|x-y|^d \int\limits_0^{\infty}
s^{-\frac{d+2}{2}}e^{-\frac{|x-y|^2}{2s}}ds
\nonumber\\ &\leq& C\phi(|x-y|) \int\limits_{0}^{\infty}
r^{\frac{d-2}{2}}e^{-r}dr\nonumber\\ &\leq& C\phi(|x-y|).
\end{eqnarray}

For $J_4(t)$, since $|x-z|>2|x-y|$, for every $\xi\in [x,y]$ (the line with endpoints $x$ and $y$)
\begin{eqnarray*}
\frac{1}{2}|x-z| \leq |\xi-z|\leq 2|x-z|.
\end{eqnarray*}
Thanks to (\ref{2.5}) and the mean value inequality, we acquire
\begin{eqnarray}\label{2.25}
|J_4(t)| &\leq&
C|x-y|\int\limits_0^t ds
\int\limits_{|x-z|> 2|x-y|}\Big(\phi(|x-z|)1_{|x-z|< r_0}+1_{|x-z|\geq r_0}|f(s,x)-f(s,z)|\Big)
\nonumber\\ &&\times (t-s)^{-\frac{(d+3)}{2}}e^{-\frac{|x-z|^2}{8(t-s)}}dz.
\end{eqnarray}
Observing that $f$ is bounded, there is a constant $C>0$ such that
\begin{eqnarray}\label{2.26}
\sup_{s\in [0,T]}|f(s,x)-f(s,z)| \leq C, \quad   {\rm for \; all} \ \ |x-z|\geq r_0.
\end{eqnarray}
Then by  (\ref{2.25}) and (\ref{2.26})
\begin{eqnarray}\label{2.27}
|J_4(t)| &\leq& C|x-y|\int\limits_{|x-z|> 2|x-y|}\Big(\phi(|x-z|)1_{|x-z|< r_0}+ 1_{|x-z|\geq r_0}\Big)|x-z|^{-d-1}dz\int\limits_0^\infty
r^{\frac{d-1}{2}}e^{-r}dr
\nonumber\\ &\leq& C|x-y|\int\limits_{|x-z|> 2|x-y|}\Big(\phi(|x-z|)1_{|x-z|< r_0}+ 1_{|x-z|\geq r_0}\Big)|x-z|^{-d-1}dz
\nonumber\\ &\leq& C|x-y|\int\limits_{2|x-y|<r< r_0}\frac{\phi(r)}{r^2}dr+C|x-y|\int\limits_{r\geq r_0}r^{-2}dr
\nonumber\\ &\leq& C|x-y|\int\limits_{2|x-y|<r<r_0}\frac{\phi(r)}{r^2}dr+C|x-y|.
\end{eqnarray}

Now combining  (\ref{2.19}), (\ref{2.20}), (\ref{2.21}) and (\ref{2.27}), for all $x,y\in\mR^d$ and $t\in [0,T]$
\begin{eqnarray}\label{2.28}
&&|\partial^2_{x_i,x_j}u(t,x)-\partial^2_{y_i,y_j}u(t,y)| \nonumber \\ &\leq &C(d,T)\Bigg[\int\limits_{r\leq 2|x-y|}\frac{\phi(r)}{r}dr+ \phi(|x-y|)+|x-y|\int\limits_{2|x-y|<r\leq r_0}\frac{\phi(r)}{r^2}dr+|x-y|\Bigg]
\nonumber \\ &\leq &C(d,T)\Bigg[\int\limits_{r\leq 2|x-y|}\frac{\phi(r)}{r}dr+ \phi(|x-y|)+2|x-y|\int\limits_{2|x-y|<r\leq r_0}\frac{\phi(r)}{r^2}dr+|x-y|\Bigg]
\nonumber \\ &= &C\Bigg[\int\limits_{r\leq |x-y|}\frac{\phi(r)}{r}dr+ \phi(|x-y|)+2|x-y|\int\limits_{|x-y|<r\leq r_0}\frac{\phi(r)}{r^2}dr+|x-y|\Bigg]\nonumber \\ &&+ \int\limits_{|x-y|<r\leq 2|x-y|}\frac{\phi(r)}{r}dr-2|x-y|\int\limits_{|x-y|<r\leq 2|x-y|}\frac{\phi(r)}{r^2}dr \nonumber \\ &\leq&C\Bigg[\int\limits_{r\leq |x-y|}\frac{\phi(r)}{r}dr+ \phi(|x-y|)+|x-y|\int\limits_{|x-y|<r\leq r_0}\frac{\phi(r)}{r^2}dr+|x-y|\Bigg],
\end{eqnarray}
which implies (\ref{2.7}).

Next we consider the case $g\neq 0$\,.   Notice that for
$g\in L^\infty([0,T];\cC_b(\mR^d;\mR^d))$ and satisfies (\ref{2.5}), $g\cdot\nabla u\in L^\infty([0,T];\cC_b(\mR^d))$ and satisfies (\ref{2.5}) with the righthand side replaced by $C\phi$. Let $\mathcal{H}$ be the set consisting of the functions in $L^\infty([0,T];\cC^2_b({\mathbb R}^d))$  with \begin{eqnarray*}\label{2.29}
|\nabla v(t,x)-\nabla v(t,y)|\leq C \phi(|x-y|), \quad {\rm for \; all }  \ x\in \mR^d, \ y\in B_{r_0}(x), \ t\in [0,T]
\end{eqnarray*}
for $v\in\mathcal{H}$ with  some constant $C>0$\,.  For $v\in\mathcal{H}$ we define a mapping
\begin{eqnarray*}
\mathcal{T}v(t,x)=\int\limits_0^{t}K(t-s,\cdot)\ast (g(s,\cdot)\cdot \nabla v(s,\cdot))(x)ds+\int\limits_0^{t}K(t-s,\cdot)\ast f(s,\cdot)(x)ds.
\end{eqnarray*}
From (\ref{2.10}),
if $f$ is in $L^\infty([0,T];\cC_b(\mR^d))$ and satisfies (\ref{2.5}), for every $x\in\mR^d$ and $y\in B_{r_0}(x)$
\begin{eqnarray}\label{2.29}
|\nabla u(t,x)-\nabla u(t,y)|&=&\Bigg|\int\limits_0^{t}\nabla K(t-s,\cdot)\ast f(s,\cdot)(x)ds-\int\limits_0^{t}\nabla K(t-s,\cdot)\ast f(s,\cdot)(y)ds\Bigg|\nonumber\\
&\leq& C\int\limits_0^{t}\int\limits_{\mR^d}\frac{K(t-s,z)}{\sqrt{t-s}}|f(s,x-z)-f(s,y-z)|dzds \nonumber\\
&\leq&C \phi(|x-y|)\sqrt{t}\leq C\phi(|x-y|).
\end{eqnarray}
Then $\mathcal{T}$ maps $\mathcal{H}$ into $\mathcal{H}$\,. Moreover  for sufficient small $T=T_{0}$, a direct verification yields that the mapping $\mathcal{T}$ is contractive.  Then there  is a unique $u\in \mathcal{H}$ such that $u=\mathcal{T}u$ and similar argument as the case $g=0$ yields   the existence and unique strong solutions of the Cauchy problem (\ref{2.1}).  Now by the classical extension technique   we construct a strong solution on $[0,T]$ for any given $T>0$ and get inequalities  (\ref{2.6})  and (\ref{2.7}) on $[0,T]$.

(ii) Let $u^n$ be the unique strong solution of (\ref{2.1}) with $f$ and $g$ replaced by  $f^n$ and $g^n$ respectively. Since $f^n$ and $g^n$ are smooth in spatial variable, $u^n$ is smooth in spatial variable as well uniformly in time.  Since   $f$ and $g$ are bounded, $f^n$ and $g^n$ are bounded with
\begin{eqnarray*}
&& \ \
\sup_n\|f^n\|_{L^\infty([0,T];\cC_b(\mR^d))}\leq \|f\|_{L^\infty([0,T];\cC_b(\mR^d))}, \\ && \sup_n\|g^n\|_{L^\infty([0,T];\cC_b(\mR^d;\mR^d))}\leq \|g\|_{L^\infty([0,T];\cC_b(\mR^d;\mR^d))}.
\end{eqnarray*}
Due to (\ref{2.6})
\begin{eqnarray*}
\sup_{n\geq 1}\|u^n\|_{L^\infty([0,T];\cC^2_b({\mathbb R}^d))} \leq C(d,T)(1+\|f\|_{L^\infty([0,T];\cC_b(\mR^d))}+
\|g\|_{L^\infty([0,T];\cC_b(\mR^d;\mR^d))}).
\end{eqnarray*}
Moreover for every $x,y\in \mR^d$ ($0<|x-y|<r_0$) and every $t\in [0,T]$
\begin{eqnarray}\label{2.30}
|f^n(t,x)-f^n(t,y)|&=&\Bigg|\int\limits_{{\mathbb R}^d}\varrho_n(z) [f(t,x-z)-f(t,y-z)]dz\Bigg|
\nonumber\\&\leq& \phi(|x-y|) \int\limits_{{\mathbb R}^d}\varrho_n(z)dz= \phi(|x-y|),
\end{eqnarray}
 and the above estimate is true for $g^n$ as well, so (\ref{2.7})  holds uniformly in $n$.

By Lemma \ref{lem2.1}
\begin{eqnarray}\label{2.31}
u_n(t,x)-u(t,x)&=&\int\limits_0^tK(t-s,\cdot)\ast [g^n(s,\cdot)\cdot \nabla
u^n(s,\cdot)-g(s,\cdot)\cdot \nabla
u(s,\cdot)](x)ds\nonumber \\ && +\int\limits_0^tK(t-s,\cdot)\ast [f^n(s,\cdot)-f(s,\cdot)](x)ds.
\end{eqnarray}
Then   by a Gronwall type argument, we arrive at
\begin{eqnarray}\label{2.32}
\| u^n-u\|_{L^\infty([0,T];\cC_b^1({\mathbb R}^d))}\leq C(\|f^n-f\|_{L^\infty([0,T];\cC_b(\mR^d))}+
\|g^n-g\|_{L^\infty([0,T];\cC_b(\mR^d;\mR^d))}),
\end{eqnarray}
where $C>0$ is independent of $n$.

On the other hand, for $f^n$ (and $g^n$) we have
\begin{eqnarray}\label{2.33}
\sup_{(t,x)\in [0,T]\times\mR^d}|f^n(t,x)-f(t,x)|&=&\sup_{(t,x)\in [0,T]\times\mR^d}\Bigg|\int\limits_{{\mathbb R}^d}\varrho_n(z) [f(t,x-z)-f(t,x)]dz\Bigg|
\nonumber\\&\leq& C \int\limits_{{\mathbb R}^d}\varrho_n(z)\phi(|z|)dz\rightarrow 0, \ as \ n\rightarrow \infty.
\end{eqnarray}
Combining (\ref{2.32}) and (\ref{2.33})
\begin{eqnarray}\label{2.34}
\lim_{n\rightarrow\infty}\| u^n-u\|_{L^\infty([0,T];\cC_b^1({\mathbb R}^d))}=0.
\end{eqnarray}

It remains to show the convergence for $\nabla^2 u^n$. For every $1\leq i,j\leq d$, every $t\in (0,T]$ and every $x\in\mR^d$
\begin{eqnarray}\label{2.35}
&&\Big|\partial^2_{x_i,x_j}u^n(t,x)-\partial^2_{x_i,x_j}u(t,x)\Big|\nonumber\\
&=&\Bigg|\int\limits_0^t\int\limits_{{\mathbb R}^d}\partial^2_{x_i,x_j}K(t-s,x-y)[g^n(s,y)\cdot\nabla u^n(s,y)-g(s,y)\cdot\nabla u(s,y)] dyds\nonumber\\
&& +\int\limits_0^t\int\limits_{{\mathbb R}^d}\partial^2_{x_i,x_j}K(t-s,x-y)[f^n(s,y)-f(s,y)]dyds\Bigg|\nonumber\\
&\leq&\Bigg|\int\limits_0^t\int\limits_{{\mathbb R}^d}\partial^2_{x_i,x_j}K(t-s,x-y)[g^n(s,y)\cdot\nabla u^n(s,y)-g(s,y)\cdot\nabla u(s,y)] dyds\Bigg|\nonumber\\
&& +\Bigg|\int\limits_0^t\int\limits_{{\mathbb R}^d}\partial^2_{x_i,x_j}K(t-s,x-y)[f^n(s,y)-f(s,y)]dyds\Bigg|
=:H^{n}_1(t,x)+H^{n}_2(t,x).
\end{eqnarray}
  First for $H^{n}_2$  we have
\begin{eqnarray}\label{2.36}
H^{n}_2(t,x)
&=& \Bigg|\int\limits_0^t\int\limits_{|x-y|\geq (t-s)^\theta }\partial^2_{x_i,x_j}K(t-s,x-y)[f^n(s,y)-f^n(s,x)-f(s,y)+f(s,x)]dyds\nonumber\\&&+\int\limits_0^t
\int\limits_{|x-y|<(t-s)^\theta}\partial^2_{x_i,x_j}K(t-s,x-y)[f^n(s,y)-f^n(s,x)-f(s,y)+f(s,x)]dyds\Bigg|
\nonumber\\
&\leq & C\|f^n-f\|_{L^\infty([0,T]\times\mR^d)}\nonumber\\&&+C\int\limits_0^t\int\limits_{|x-y|<(t-s)^\theta\wedge r_0}\frac{1}{t-s}K(t-s,x-y)|f^n(s,y)-f^n(s,x)-f(s,y)+f(s,x)|dyds\nonumber\\&&+C\int\limits_0^t
\int\limits_{(t-s)^\theta\wedge r_0 \leq |x-y|\leq (t-s)^\theta
}\frac{1}{t-s}K(t-s,x-y)\nonumber\\ &&\quad \times |f^n(s,y)-f^n(s,x)-f(s,y)+f(s,x)|dyds,
\end{eqnarray}
for  $\theta\in (0,1/2)$.
By (\ref{2.5}) and (\ref{2.30}),
\begin{eqnarray*}
|f^n(s,y)-f^n(s,x)-f(s,y)+f(s,x)|\leq C\phi(|x-y|), \quad  {\rm for \; all}\ \  |x-y|<r_0.
\end{eqnarray*}
From (\ref{2.16})
\begin{eqnarray*}
\frac{1}{t-\cdot}K(t-\cdot,\cdot-y) \phi(|\cdot-x|) \in L^1(O), \quad  O=\{(s,y)| \ s\in [0,t], \ |x-y|< (t-s)^\alpha\wedge r_0\}.
\end{eqnarray*}
By the assumption on $f$ and the definition of $f_n$, the integrand of the last term in (\ref{2.36}) is integrable.  Then by the dominated convergence theorem and (\ref{2.33}), we have
\begin{eqnarray*}
\lim_{n\rightarrow\infty}\sup_{(t,x)\in [0,T]\times\mR^d}H^{n}_{2}(t, x)=0.
\end{eqnarray*}
Noticing  that $g^n(s,y)\cdot\nabla u^n(s,y)$ and $g(s,y)\cdot\nabla u(s,y)$ are of class $L^\infty([0,T];\cC_b(\mR^d))$, satisfy (\ref{2.5}) and
\begin{eqnarray*}
\lim_{n\rightarrow\infty}\|g^n\cdot\nabla u^n-g\cdot\nabla u\|_{L^\infty([0,T];\cC_b({\mathbb R}^d))}=0.
\end{eqnarray*}
So the argument of convergence for $H^{n}_1$ is  same as  that  for $H^{n}_2$, then we have
\begin{eqnarray}\label{2.37}
\lim_{n\rightarrow\infty}\|\nabla^2u^n-\nabla^2u\|_{L^\infty([0,T];\cC_b({\mathbb R}^d))}=0,
\end{eqnarray}
and combining (\ref{2.34}) and (\ref{2.37}), we show (\ref{2.9}). $\Box$

\begin{remark} \label{rem2.1}
There are some further properties of  $\nabla^2 u$ and $\partial_tu$  by some estimates on the righthand side of (\ref{2.7}). We consider   a nonnegative Dini  function $\phi$.  For every $r_0>0$, then
\begin{eqnarray}\label{2.38}
\int\limits_0^{r_0}\frac{\phi(r)}{r}dr<+\infty.
\end{eqnarray}
On the other hand,
\begin{eqnarray*}
\varepsilon\int\limits_{\varepsilon<r\leq r_0}\frac{\phi(r)}{r^2}dr=\int\limits_{0<r\leq r_0}\frac{\phi(r)}{r^2}1_{(\varepsilon,r_0)}(r)\varepsilon dr=:\int\limits_{0<r\leq r_0}h_\varepsilon(r) dr.
\end{eqnarray*}
Since $h_\varepsilon(r)\leq \phi(r)/r$, by the dominated convergence theorem,
\begin{eqnarray}\label{2.39}
\lim_{\varepsilon\rightarrow 0}\varepsilon\int\limits_{\varepsilon<r\leq r_0}\frac{\phi(r)}{r^2}dr=\int\limits_{0<r\leq r_0}\lim_{\varepsilon\rightarrow 0}h_\varepsilon(r) dr=0.
\end{eqnarray}
By (\ref{2.39}), all the terms in the righthand side of (\ref{2.7}) are meaningful. Moreover, by the above estimates,  $\nabla^2 u$ and $\partial_tu$ are uniformly continuous in spatial variable uniformly in time.
\end{remark}

\begin{remark} \label{rem2.2}
Theorem \ref{the2.1} generalizes the existing results for  the Cauchy problem (\ref{2.1}).
By the classical parabolic theory (\cite[Ch. 4]{LSU}, \cite{Kry}), if $f\in L^p([0,T]\times{\mathbb R}^d),\;  g\in L^p([0,T]\times{\mathbb R}^d;\mR^d)$ with $p\in (1,+\infty)$,  $u\in L^p([0,T];W^{2,p}({\mathbb R}^d))\cap W^{1,p}([0,T];L^p({\mathbb R}^d))$. However, if $f\in L^\infty([0,T]\times{\mathbb R}^d), g\in L^\infty([0,T]\times{\mathbb R}^d;\mR^d)$, in general $u\not\in L^\infty([0,T];W^{2,\infty}({\mathbb R}^d))\cap W^{1,\infty}([0,T];L^\infty({\mathbb R}^d))$.
Nevertheless,  recent result~\cite[Lemma 2.1]{WDGL} shows that if $f\in L^p([0,T];{\mathcal C}^\varsigma_b({\mathbb R}^d)), g\in L^p([0,T];{\mathcal C}^\varsigma_b({\mathbb R}^d;\mR^d))$ with $p\in (1,+\infty], \varsigma\in (0,1)$, and $\varsigma>2/p$
\begin{eqnarray*}
u\in L^\infty([0,T];{\mathcal C}^{2,\varsigma-\frac{2}{p}}_b({\mathbb R}^d))\cap W^{1,\infty}([0,T];{\mathcal C}^{\varsigma-\frac{2}{p}}_b({\mathbb R}^d)).
\end{eqnarray*}
Here, by assuming the Dini continuity on $f$ and $g$, we establish the $W^{2,\infty}$ estimates as well. Moreover, the second derivatives for spatial variable are also uniformly continuous.
\end{remark}

\begin{remark} \label{rem2.3}
(i) From (\ref{2.7}), if there exists $\eta>0$ such that for $r>0$ small enough
\begin{eqnarray}\label{2.40}
\eta\int\limits_0^r\frac{\phi(s)}{s}ds=\phi(r)
\end{eqnarray}
and
\begin{eqnarray}\label{2.41}
r\int\limits_{r<s\leq r_0}\frac{\phi(s)}{s^2}ds\leq C\phi(r),
\end{eqnarray}
then the maximum Dini regularity for (\ref{2.1}) holds. Now, we solve the integral equation (\ref{2.40}) and get $\phi(r)=C_0r^\eta$. Since the Dini function $\phi$ is given in (\ref{2.5}), $\eta\in (0,1)$ is the best choice. Then,
\begin{eqnarray*}
r\int\limits_{r<s\leq r_0}\frac{\phi(s)}{s^2}ds=\frac{C_0}{1-\eta}[ r^\eta-rr_0^{\eta-1}]\leq \frac{\phi(r)}{1-\eta},
\end{eqnarray*}
for $r$ is sufficiently small. And now, we recover the classical Schauder theory for parabolic equations of second order.

(ii) If $\phi$ is only Dini continuous but not H\"{o}lder continuous, the functions given in the first and third terms on the right hand side of (\ref{2.7}) will not preserve the same regularity as $\phi$, and thus   the maximum regularity for (\ref{2.1}) may be not true. For example, we choose $r_0=1/2$ and $\phi(r)=C|\log(r)|^{-\alpha}$ with some $\alpha>1$, then there exist two positive real numbers $C_1(d,\alpha)$ and $C_2(d,\alpha)$ such that for $|x-y|>0$ sufficiently small
\begin{eqnarray*}
\frac{C_1(d,\alpha)}{|\log(|x-y|)|^{\alpha-1}}\leq \int\limits_{r\leq |x-y|}\frac{\phi(r)}{r}dr+|x-y|\int\limits_{|x-y|<r\leq r_0}\frac{\phi(r)}{r^2}dr\leq \frac{C_2(d,\alpha)}{|\log(|x-y|)|^{\alpha-1}}.
\end{eqnarray*}
However, $\phi(r)=o(|\log(r)|^{1-\alpha})$ as $r\rightarrow 0$. Therefore, if $\alpha\in (1,2)$, $\nabla^2u(t,\cdot)$ may be only uniformly continuous but not Dini continuous.
\end{remark}

To make the maximum regularity for (\ref{2.1}) true, we need introduce the following notion.
\begin{definition} \label{def2.2} Let $\psi: \mR_+\rightarrow \bar{{\mathbb R}}_+$ be a monotone continuous function. Suppose $\phi: \mR_+\rightarrow \mR_+$ is continuous and $\phi(0)=0$.

(i) We call $\phi$ a H\"{o}lder-Dini function with H\"{o}lder exponent $\vartheta\in (0,1)$ if $r^{-\vartheta}\phi(r)=\psi(r)$ is Dini but not H\"{o}lder continuous on $\mR_+$;

(ii) We call $\phi$ a strong H\"{o}lder function with H\"{o}lder exponent $\vartheta\in (0,1)$ if

(1) for every $\zeta\in (0,\vartheta)$, $r^{-\zeta}\phi(r)$ is H\"{o}lder continuous on $\mR_+$;

(2) $r^{-\vartheta}\phi(r)=\psi(r)$ is strictly monotone increasing and continuous but not H\"{o}lder continuous on $\mR_+$;

(3) $\psi(r)\downarrow0$ as $r \downarrow 0$;

(iii) We call $\phi$ a weak H\"{o}lder function with H\"{o}lder exponent $\vartheta\in (0,1)$ if

(1) for every $\zeta\in (0,\vartheta)$, $r^{-\zeta}\phi(r)$ is H\"{o}lder continuous on $\mR_+$ and $r^{-\zeta}\phi(r)\rightarrow 0$ as $r \downarrow 0$;

(2) $r^{-\vartheta}\phi(r)=\psi(r)$ is strictly monotone decreasing and continuous for $r\in (0,\infty)$;

(3) $\psi(r)\uparrow +\infty$ as $r \downarrow 0$.

A measurable function $h:\mR^{d}\rightarrow \mR$  is said to be H\"{o}lder-Dini (strong H\"{o}lder or weak H\"{o}lder) continuous with H\"{o}lder exponent $\vartheta\in (0,1)$ if
\begin{eqnarray}\label{2.42}
|h(x)-h(y)|\leq \phi(|x-y|),
\end{eqnarray}
where $\phi$ is a H\"{o}lder-Dini (strong H\"{o}lder or weak H\"{o}lder) function with H\"{o}lder exponent $\vartheta$.
\end{definition}

If $f$ and $g$ are H\"{o}lder-Dini or strong H\"{o}lder or weak H\"{o}lder continuous, then the maximum  regularity for (\ref{2.1}) is still true.
\begin{corollary} \label{cor2.1} Let $\phi$ be a H\"{o}lder-Dini or strong H\"{o}lder or weak H\"{o}lder function with H\"{o}lder exponent $\vartheta\in (0,1)$. Let $f\in L^\infty([0,T];\cC_b(\mR^d))$ and $g\in L^\infty([0,T];\cC_b(\mR^d;\mR^d))$ such that (\ref{2.5}) holds. Then there exists a unique strong solution $u$ to (\ref{2.1}). Moreover, if $\phi^\prime$ is also continuous in $(0,\infty)$, then for every $1\leq i,j\leq d$ and every $x,y\in \mR^d$, there is a constant $C(d,T)$ such that
\begin{eqnarray}\label{2.43}
|\partial^2_{x_i,x_j}u(t,x)-\partial^2_{y_i,y_j}u(t,y)|\leq C(d,T)\phi(|x-y|), \quad    for \;  all  \ t\in [0,T].
\end{eqnarray}
\end{corollary}

Before proving the above result, we need a useful lemma.

\begin{lemma} [L'Hospital's rule \cite{Daw}, p. 346] \label{lem2.2} Suppose that we have one of the following cases,
\begin{eqnarray*}
\lim_{r\rightarrow \lambda}\frac{f_1(r)}{f_2(r)} = \frac{0}{0} \quad or \quad  \lim_{r\rightarrow \lambda}\frac{f_1(r)}{f_2(r)} = \frac{\pm \infty}{\pm \infty} ,
\end{eqnarray*}
where $\lambda$ can be any real number, infinity or negative infinity, In these cases, we have
\begin{eqnarray*}
\lim_{r\rightarrow \lambda}\frac{f_1(r)}{f_2(r)} =\lim_{r\rightarrow \lambda}\frac{f_1^\prime(r)}{f_2^\prime(r)}.
\end{eqnarray*}
\end{lemma}

\medskip\noindent\textbf{Proof of Corollary \ref{cor2.1}.} We only need to prove (\ref{2.43}) for $|x-y|$ sufficiently small. By (\ref{2.7}) for every $1\leq i,j\leq d$ and every $x,y\in \mR^d$ ($|x-y|<r_0$), there is a constant $C(d,T)$ such that
\begin{eqnarray}\label{2.44}
&&|\partial^2_{x_i,x_j}u(t,x)-\partial^2_{y_i,y_j}u(t,y)| \nonumber \\ &\leq &C\left[\int\limits_{r\leq |x-y|}\frac{\phi(r)}{r}dr+ \phi(|x-y|)+|x-y|\int\limits_{|x-y|<r\leq r_0}\frac{\phi(r)}{r^2}dr\right].
\end{eqnarray}
From (\ref{2.44}) if
\begin{eqnarray}\label{2.45}
\limsup_{r\rightarrow 0} \frac{\int\limits_0^r\frac{\phi(s)}{s}ds}{\phi(r)}<+\infty\quad
\mbox{and} \quad
\limsup_{r\rightarrow 0} \frac{r\int\limits_r^{r_0}\frac{\phi(s)}{s^2}ds}{\phi(r)} <+\infty,
\end{eqnarray}
then for $|x-y|$ small enough
\begin{eqnarray*}
\max\{|x-y|\int\limits_{|x-y|<r\leq r_0}\frac{\phi(r)}{r^2}dr, \int\limits_{r\leq |x-y|}\frac{\phi(r)}{r}dr\}\leq C \phi(|x-y|),
\end{eqnarray*}
thus (\ref{2.43}) is proved.

Now let us check (\ref{2.45}). Since  $\phi$ is a H\"{o}lder-Dini or strong H\"{o}lder or weak H\"{o}lder function with H\"{o}lder exponent $\vartheta\in (0,1)$, we can rewrite $\phi(r)$ by $r^{\vartheta}\psi(r)$ with $\psi$ a monotone continuous function. Applying Lemma \ref{lem2.2}, then
\begin{eqnarray}\label{2.46}
\lim_{r\rightarrow 0} \frac{\int\limits_0^r\frac{\phi(s)}{s}ds}{\phi(r)}=\lim_{r\rightarrow 0} \frac{\int\limits_0^rs^{\vartheta-1}\psi(s)ds}{r^{\vartheta}\psi(r)}=\lim_{r\rightarrow 0} \Big[\vartheta+\frac{r\psi^\prime(r)}{\psi(r)}\Big]^{-1}
\end{eqnarray}
and
\begin{eqnarray}\label{2.47}
\lim_{r\rightarrow 0} \frac{r\int\limits_r^{r_0}s^{\vartheta-2}\psi(s)ds}{r^{\vartheta}\psi(r)}=\lim_{r\rightarrow 0} \frac{\int\limits_r^{r_0}s^{\vartheta-2}\psi(s)ds}{r^{\vartheta-1}\psi(r)}=\lim_{r\rightarrow 0} \Big[1-\vartheta-\frac{r\psi^\prime(r)}{\psi(r)}\Big]^{-1}.
\end{eqnarray}

If $\phi$ is a H\"{o}lder-Dini or strong H\"{o}lder function, then
\begin{eqnarray}\label{2.48}
\lim_{r\rightarrow 0} \log(\psi(r))=-\infty,
\end{eqnarray}
and if $\phi$ is a weak H\"{o}lder function, then
\begin{eqnarray}\label{2.49}
\lim_{r\rightarrow 0} \log(\psi(r))=+\infty.
\end{eqnarray}
From (\ref{2.48}) and (\ref{2.49}), by using Lemma \ref{lem2.2} again, we gain
\begin{eqnarray}\label{2.50}
\lim_{r\rightarrow 0} \frac{r\psi^\prime(r)}{\psi(r)}=
\left\{\begin{array}{ll}
\lim\limits_{r\rightarrow 0} \frac{\log(\psi(r))}{\log(r)}, & \hbox{when $\phi$ is H\"{o}lder-Dini or strong H\"{o}lder;} \\
-\lim\limits_{r\rightarrow 0} \frac{\log(\psi(r))}{\log(r^{-1})}, & \hbox{when $\phi$ is weak H\"{o}lder.}\end{array}\right.
\end{eqnarray}

When $\phi$ is a H\"{o}lder-Dini or strong H\"{o}lder function, then for every $\zeta>0$, $\psi(r)\geq r^{\zeta}$ as $r\rightarrow 0$.
Thus, there exists $\delta=\delta(\zeta)>0$ such that
\begin{eqnarray*}
[\psi(r)]^{-1}\leq r^{-\zeta}, \quad  for\; all  \ r\in (0,\delta].
\end{eqnarray*}
It follows that
\begin{eqnarray}\label{2.51}
\limsup_{r\rightarrow 0} \frac{\log(\psi(r))}{\log(r)}=\limsup_{r\rightarrow 0} \frac{\log([\psi(r)]^{-1})}{\log(r^{-1})}\leq \limsup_{r\rightarrow 0} \frac{\log(r^{-\zeta})}{\log(r^{-1})}=\zeta.
\end{eqnarray}
Since $\zeta>0$ is arbitrary and $\log(\psi(r))/\log(r)\geq 0$, from (\ref{2.50}) and (\ref{2.51}), we conclude
\begin{eqnarray}\label{2.53}
\lim_{r\rightarrow 0} \frac{r\psi^\prime(r)}{\psi(r)}=\lim_{r\rightarrow 0} \frac{\log(\psi(r))}{\log(r)}=0.
\end{eqnarray}

When $\phi$ is a weak H\"{o}lder function, then for every $\zeta_1>0$, $r^{\zeta_1}\psi(r)\rightarrow0$ as $r\rightarrow 0$. Thus, there exists $\delta_1=\delta_1(\zeta_1)>0$ such that
\begin{eqnarray*}
\psi(r)\leq r^{-\zeta_1}, \quad  for\; all  \ r\in (0,\delta_1],
\end{eqnarray*}
which implies that
\begin{eqnarray}\label{2.54}
\limsup_{r\rightarrow 0} \frac{\log(\psi(r))}{\log(r^{-1})}\leq \limsup_{r\rightarrow 0} \frac{\log(r^{-\zeta_1})}{\log(r^{-1})}=\zeta_1.
\end{eqnarray}
Since $\zeta_1>0$ is arbitrary and $\log(\psi(r))/\log(r^{-1})\geq 0$, from (\ref{2.50}) and (\ref{2.54}), we conclude
\begin{eqnarray}\label{2.55}
\lim_{r\rightarrow 0} \frac{r\psi^\prime(r)}{\psi(r)}=-\lim_{r\rightarrow 0} \frac{\log(\psi(r))}{\log(r^{-1})}=0.
\end{eqnarray}

Combining (\ref{2.53}) and (\ref{2.55}), then (\ref{2.45}) holds. From this we complete the proof. $\Box$

\begin{example} \label{exa2.1} We choose $\phi(r)=Cr^\vartheta|\log(r)|^\alpha$ with $\vartheta\in (0,1)$ and $\alpha\in {\mathbb R}$, then $\phi$ is a H\"{o}lder-Dini function if $\alpha<-1$, a strong H\"{o}lder function if $\alpha<0$, H\"{o}lder continuous if $\alpha=0$, a weak H\"{o}lder function if $\alpha>0$. Using Lemma \ref{lem2.2}, we have
\begin{eqnarray*}
\lim_{r\rightarrow 0} \frac{\int\limits_0^r\frac{\phi(s)}{s}ds}{\phi(r)}=\lim_{r\rightarrow 0} \frac{\phi(r)}{r\phi^\prime(r)}=\lim_{r\rightarrow 0} \frac{\phi(r)}{\phi(r)[\vartheta-\alpha |\log(r)|^{-1}]}
=\frac{1}{\vartheta}
\end{eqnarray*}
and
\begin{eqnarray*}
\lim_{r\rightarrow 0} \frac{r\int\limits_r^{r_0}\frac{\phi(s)}{s^2}ds}{\phi(r)}&=&-\lim_{r\rightarrow 0} \frac{\phi(r)}{r\phi^\prime(r)}+ \lim_{r\rightarrow 0} \frac{\int\limits_r^{r_0}\frac{\phi(s)}{s^2}ds}{\phi^\prime(r)}
\nonumber\\ &=&-\frac{1}{\vartheta}-\lim_{r\rightarrow 0} \frac{\phi(r)}{r^2\phi^{\prime\prime}(r)}
\nonumber\\ &=&-\frac{1}{\vartheta}-\lim_{r\rightarrow 0} \frac{\phi(r)}{\phi(r)[\vartheta(\vartheta-1)-\alpha (2\vartheta-1)|\log(r)|^{-1}-\alpha(1-\alpha)|\log(r)|^{-2}]}\nonumber\\ &=&-\frac{1}{\vartheta}-\frac{1}{\vartheta(\vartheta-1)}=\frac{1}{1-\vartheta}.
\end{eqnarray*}
Therefore, for $r$ small enough, we gain
\begin{eqnarray}\label{2.56}
\int\limits_0^r\frac{\phi(s)}{s}ds\leq \frac{1+\vartheta}{\vartheta} \phi(r), \quad
r\int\limits_r^{r_0}\frac{\phi(s)}{s^2}ds\leq
\frac{2-\vartheta}{1-\vartheta}\phi(r).
\end{eqnarray}
Combining (\ref{2.7}) and (\ref{2.56}), the maximum regularity for (\ref{2.1}) is preserved. This result as we know is new.
\end{example}

\section{Stochastic flows for SDEs with  bounded and Dini drift}\label{sec3}\setcounter{equation}{0}

Given real number  $T>0$\,, for $s\in[0,T]$ and $x\in {\mathbb R}^d$, consider the  following SDE
\begin{eqnarray}\label{3.1}
dX(s,t)=b(t,X(s,t))dt+dB(t), \quad  t\in(s,T], \quad  X(s,t)|_{t=s}=x\,.
\end{eqnarray}
We intend to show the existence of a stochastic flow for equation (\ref{3.1})\,. First we give the following definition.

\begin{definition} [\cite{Kun2}, p. 114] \label{def2.1}
A stochastic homeomorphisms flow of class $\mathcal{C}^{\beta}$ with $\beta\in (0,1)$ on
$(\Omega, \mathcal{F},{\mathbb P}, (\mathcal{F}_t)_{0\leq t\leq T})$ associated to (\ref{3.1}) is a
map $(s,t,x,\omega) \rightarrow X(s,t,x)(\omega)$, defined for
$0\leq s \leq t \leq T, \ x\in {\mathbb R}^d, \ \omega \in \Omega$ with
values in ${\mathbb R}^d$, such that

(i)  the process
$\{X(s,\cdot,x)\}= \{X(s,t,x), \ t\in [s,T]\}$ is a continuous
$\{\mathcal{F}_{s,t}\}_{s\leq t\leq T}$-adapted solution of (\ref{3.1}),  for  every   $s\in [0,T],\  x \in {\mathbb R}^d$;

(ii)  the functions $X(s,t,x)$ and $X^{-1}(s,t,x)$ are continuous in $(s,t,x)$ and are of class ${\mathcal C}^\beta$ in $x$ uniformly in $(s,t)$\,, ${\mathbb P}$-a.s., for all  $0\leq s \leq t \leq T$;

(iii)   $X(s,t,x)=X(r,t,X(s,r,x))$  for all
$0\leq s\leq r \leq t \leq T$, $x\in {\mathbb R}^d$\,,  ${\mathbb P}$-a.s.,  and $X(s,s,x)=x$\,.

Further if

(iv)   for all  $0\leq s \leq t \leq T$, the functions
$\nabla X(s,t,x)$ and $\nabla X^{-1}(s,t,x)$ are continuous in $(s,t,x)$\,, ${\mathbb P}$-a.s.,

  it is called a stochastic diffeomorphisms flow.
\end{definition}

Now we  state the  main result of the section.
\begin{theorem} \label{the3.1} Let $b\in L^\infty([0,T];\cC_b(\mR^d;\mR^d))$ with integer  $d\geq 1$. Suppose that $r_0\in (0,1)$ and there is a Dini function $\phi$ such that for every $x\in\mR^d$
\begin{eqnarray}\label{3.2}
|b(t,x)-b(t,y)|\leq \phi(|x-y|), \quad  {\rm for\; all}  \ \ y\in B_{r_0}(x), \ t\in [0,T].
\end{eqnarray}
In addition that for every $p\geq 1$\,,  there is a small enough positive real number $\delta=\delta(p)<r_0$ such that $F_\delta^p(\cdot)$ is increasing and concave on $[0,\delta]$, where
\begin{eqnarray}\label{3.3}
F_\delta(r)=\int\limits_{s\leq r}\frac{\phi(s)}{s}ds+ \phi(r)+r\int\limits_{r<s\leq \delta}\frac{\phi(s)}{s^2}ds+r, \quad    \ r\in [0,\delta].
\end{eqnarray}

(i) For
every $s\in [0,T ]$ and $ x\in {\mathbb R}^d$, SDE
(\ref{3.1}) has a unique continuous adapted solution $\{X(s,t,x)(\omega), \
t\in [s,T ], \ \omega \in\Omega\}$, which forms a stochastic diffeomorphisms flow. For every $p\geq 1$, there is a constant $C(p,d,T)>0$  such that
\begin{eqnarray}\label{3.4}
\sup_{0\leq s\leq
T}{\mathbb E}\sup_{s\leq t \leq T}|X(s,t,x)|^p+\sup_{x\in\mR^d}\sup_{0\leq s\leq
T}{\mathbb E}\sup_{s\leq t \leq T}\|\nabla X(s,t,x)\|^p\leq C(p,d,T).
\end{eqnarray}
Moreover, for every $p\geq 1$ and every $x,y\in \mR^d$,
\begin{eqnarray}\label{3.5}
\sup_{0\leq s\leq
T}{\mathbb E}[\sup_{s\leq t \leq T}|X(s,t,x)-X(s,t,y)|^p]\leq C|x-y|^p
\end{eqnarray}
and
\begin{eqnarray}\label{3.6}
&&\sup_{0\leq s\leq
T}{\mathbb E}[\sup_{s\leq t \leq T}\|\nabla X(s,t,x)-\nabla X(s,t,y)\|^p]
\nonumber \\
&\leq &C\Bigg[\int\limits_{r\leq |x-y|}\frac{\phi(r)}{r}dr+ \phi(|x-y|)+|x-y|\int\limits_{|x-y|<r\leq r_0}\frac{\phi(r)}{r^2}dr\Bigg]^p1_{|x-y|<r_0}+C|x-y|^p.
\end{eqnarray}

(ii) Let $\varrho_n$ be given in (\ref{2.8}) and  $b^n(t,x)=b\ast \varrho_n(t,x)$. Let $X^n$ be the  stochastic flows corresponding to the vector field $b^n$. Then for every $p\geq 1$,
\begin{eqnarray}\label{3.7}
\lim _{n\rightarrow\infty}\sup_{x\in {\mathbb R}^d} \sup_{0\leq s\leq
T}{\mathbb E}\left[\sup_{s\leq t \leq T}|X^n(s,t,x)-X(s,t,x)|^p\right]=0
\end{eqnarray}
and
\begin{eqnarray}\label{3.8}
\lim_{n\rightarrow \infty}\sup_{x\in {\mathbb R}^d} \sup_{0\leq s\leq
T}{\mathbb E}\left[\sup_{s\leq t \leq T}\|\nabla X^n(s,t,x)-\nabla X(s,t,x)\|^p\right]=0.
\end{eqnarray}
\end{theorem}
\noindent
\textbf{Proof.} (i) Let $\lambda>0$ be a real number. Consider the following backward heat equation
\begin{eqnarray}\label{3.9}
\left\{
\begin{array}{ll}
\partial_{t}U(t,x) +\frac{1}{2}\Delta U(t,x)+b(t,x)\cdot \nabla U(t,x)\\ \quad \quad=
\lambda U(t,x)-b(t,x), \quad  (t,x)\in [0,T)\times {\mathbb R}^d, \\
U(T,x)=0, \quad   x\in{\mathbb R}^d.
  \end{array}
\right.
\end{eqnarray}
By  Theorem \ref{the2.1}, the above Cauchy problem has  a unique solution $U\in L^\infty([0,T];\cC^2_b({\mathbb R}^d;{\mathbb R}^d))\cap W^{1,\infty}([0,T];\mathcal{C}_b({\mathbb R}^d;{\mathbb R}^d))$.
Moreover, (\ref{2.6}) is true and the constant $C$ in (\ref{2.6}) is independent of $\lambda$. Further, for every $x,y\in \mR^d$ and $t\in [0,T)$, there is a constant $C(d,T)$ such that
\begin{eqnarray}\label{3.10}
&&|\partial^2_{x_i,x_j}U(t,x)-\partial^2_{y_i,y_j}U(t,y)| \nonumber \\ &\leq &C(d,T)\Bigg[\int\limits_{r\leq |x-y|}\frac{\phi(r)}{r}dr+ \phi(|x-y|)+|x-y|\int\limits_{|x-y|<r\leq r_0}\frac{\phi(r)}{r^2}dr\Bigg]1_{|x-y|<r_0}\nonumber \\ &&+C(d,T)|x-y|.
\end{eqnarray}

By Lemma \ref{lem2.1}, the unique strong solution $U$ has the following  representation
\begin{eqnarray*}
U(t,x)=\int\limits_0^{T-t}e^{-\lambda r}K(r,\cdot)\ast[ b(t+r,\cdot)+b(t+r,\cdot)\cdot \nabla U(t+r,\cdot)](x)dr.
\end{eqnarray*}
For every $1\leq i\leq d$, every $x\in \mR^d$ and $t\in [0,T]$
\begin{eqnarray}\label{3.11}
|\partial_{x_i}U(t,x)|&=&\Bigg|\int\limits_0^{T-t}e^{-\lambda r}dr\int\limits_{{\mathbb R}^d}\partial_{x_i}K(r,x-z)[b(t+\tau,z)+b(t+r,z)\cdot \nabla U(t+r,z)]dz\Bigg|
\nonumber\\&\leq&C\|b\|_{L^\infty([0,T];\mathcal{C}_b({\mathbb R}^d;{\mathbb R}^d))}(1+\|b\|_{L^\infty([0,T];\mathcal{C}_b({\mathbb R}^d;{\mathbb R}^d))})
\int\limits_0^Tr^{-\frac{1}{2}}e^{-\lambda r} dr
\nonumber\\&\leq&C\|b\|_{L^\infty([0,T];\mathcal{C}_b({\mathbb R}^d;{\mathbb R}^d))}(1+\|b\|_{L^\infty([0,T];\mathcal{C}_b({\mathbb R}^d;{\mathbb R}^d))})\lambda^{-\frac{1}{3}},
\end{eqnarray}
where the H\"{o}lder inequality is applied  in the last inequality.
Then  letting  $\lambda$ tend to infinity in~(\ref{3.11}) yields
\begin{eqnarray}\label{3.12}
\|\nabla U\|_{L^\infty([0,T];{\mathcal C}_b({\mathbb R}^d;\mR^{d\times d}))}\longrightarrow 0.
\end{eqnarray}
Therefore, there is a large real number $\lambda_0>0$ such that if $\lambda\geq \lambda_0$
\begin{eqnarray}\label{3.13}
\|\nabla U\|_{L^\infty([0,T];{\mathcal C}_b({\mathbb R}^d;\mR^{d\times d}))}\leq \frac12.
\end{eqnarray}
Now  for  $\lambda\geq \lambda_{0}$,  define $\gamma(t,x)=x+U(t,x)$
\begin{eqnarray*}
\frac12\leq \|\nabla \gamma\|_{L^\infty([0,T];{\mathcal C}_b({\mathbb R}^d;\mR^{d\times d}))}\leq \frac32.
\end{eqnarray*}
Then  $\gamma(t)$ forms a nonsingular diffeomorphism of class $\mathcal{C}^2$ uniformly in $t\in [0,T]$ by the classical Hadamard theorem (\cite[p.330]{Protter}).  Moreover, for every $t\in [0,T]$, the inverse of $\gamma(t)$ (denoted by~$\gamma^{-1}(t)$) has bounded first and second spatial derivatives, uniformly in $t\in [0,T]$, and
\begin{eqnarray}\label{3.14}
\frac23\leq\|\nabla \gamma^{-1}\|_{L^\infty([0,T];{\mathcal C}_b({\mathbb R}^d;\mR^{d\times d}))}\leq 2.
\end{eqnarray}

Noticing that  if $X(s,t)$ satisfies SDE (\ref{3.1}),   $Y(s,t)=X(s,t)+U(t,X(s,t))$ $(=:\gamma(t,X(s,t)))$ satisfies the following SDE
\begin{eqnarray}\label{3.15}
\left\{
  \begin{array}{ll}
d Y(s,t)=\lambda U(t,\gamma^{-1}(t,Y(s,t)))dt+ [I+\nabla
U(t,\gamma^{-1}(t,Y(s,t)))] dB(t)\\ \quad\quad\quad\ \ =:\tilde{b}(t,Y(s,t))dt+\tilde{\sigma}(t,Y(s,t))dB(t), \ t\in(s,T], \\
Y(s,t)|_{t=s}=y,
  \end{array}
\right.
\end{eqnarray}
and vice versa, (\ref{3.15}) and  (\ref{3.1}) are equivalent. Now for  equation (\ref{3.15})\,, the coefficients are globally Lipschitz continuous,   by the result of Kunita~\cite[Theorem 4.3, p.227]{Kun1}\,,  there is  a unique stochastic homeomorphisms flow of class $\cC^\beta \ (\beta\in (0,1))$ defined by $Y(s,t, y)$\,.  Moreover, for every $x,y\in \mR^d$ and every $p\geq 1$, there is a constant $C(p,T)>0$ such that
\begin{eqnarray}\label{3.16}
\sup_{0\leq s\leq
T}{\mathbb E}\sup_{s\leq t \leq T}|Y(s,t,y)|^p\leq C(p,T)
\end{eqnarray}
and
\begin{eqnarray}\label{3.17}
\sup_{0\leq s\leq
T}{\mathbb E}[\sup_{s\leq t \leq T}|Y(s,t,x)-Y(s,t,y)|^p]\leq C(p,T)|x-y|^p.
\end{eqnarray}
On the other hand, by   $X(s,t)=\gamma^{-1}(t,Y(s,t))$\,,  (\ref{3.1}) also defines a unique stochastic homeomorphisms flow of class $\cC^\beta \ (\beta\in (0,1))$.  Moreover by (\ref{3.16}) and (\ref{3.17}) for every $p\geq 1$ and every $x,y\in\mR^d$
\begin{eqnarray}\label{3.18}
&&\sup_{0\leq s\leq
T}{\mathbb E}\sup_{s\leq t \leq T}|X(s,t,x)|^p\nonumber \\ &\leq& 2^{p-1}\sup_{0\leq s\leq
T}{\mathbb E}\sup_{s\leq t \leq T}[|Y(s,t,\gamma(s,x))|^p+|\gamma^{-1}(t,Y(s,t,\gamma(s,x)))-\gamma^{-1}(t,\gamma(s,x))|^p]
\nonumber \\ &\leq& C(p,T)+2^{p-1}\|\nabla \gamma^{-1}\|^p_{L^\infty([0,T];\cC_b(\mR^d;\mR^{d\times d}))} \sup_{0\leq s\leq
T}{\mathbb E}\sup_{s\leq t \leq T}|Y(s,t,\gamma(s,x))-\gamma(s,x)|^p\nonumber \\ &\leq& C(p,d,T)
\end{eqnarray}
and
\begin{eqnarray}\label{3.19}
&&\sup_{0\leq s\leq
T}{\mathbb E}[\sup_{s\leq t \leq T}|X(s,t,x)-X(s,t,y)|^p]\nonumber \\ &=& \sup_{0\leq s\leq
T}{\mathbb E}[\sup_{s\leq t \leq T}|\gamma^{-1}(t,Y(s,t,\gamma(s,x))) -\gamma^{-1}(t,Y(s,t,\gamma(s,y)))|^p]
\nonumber \\ &\leq& C \|\nabla \gamma^{-1}\|^p_{L^\infty([0,T];\cC_b(\mR^d;\mR^{d\times d}))}  \sup_{0\leq s\leq
T}{\mathbb E}[\sup_{s\leq t \leq T}|Y(s,t,\gamma(s,x)) -Y(s,t,\gamma(s,y))|^p]
\nonumber \\ &\leq& C \|\nabla \gamma^{-1}\|^p_{L^\infty([0,T];\cC_b(\mR^d;\mR^{d\times d}))}  |\gamma(s,x) -\gamma(s,y)|^p\nonumber \\ &\leq& C \|\nabla \gamma^{-1}\|^p_{L^\infty([0,T];\cC_b(\mR^d;\mR^{d\times d}))}\|\nabla \gamma\|^p_{L^\infty([0,T];\cC_b(\mR^d;\mR^{d\times d}))}   |x -y|^p\nonumber \\ &\leq&
C(p,d,T)|x-y|^p.
\end{eqnarray}

Now we just need to prove the continuity of $\nabla_x X(s,t,x)$ and the inequality (\ref{3.6}). For this  we  consider $\nabla_y Y(s,t,y)$\,.  First the stochastic flow $\{X(s,t)(\cdot)\}$ is weak differentiable, that is,~$X(s,t)(\cdot)\in L^2(\Omega;W^{1,p}_{loc}(\mR^d;\mR^d))$~(\cite[Theorem 3]{MNP}), so does  the stochastic flow $\{Y(s,t)(\cdot)\}$.  Then differentiate~$Y(s,t)$ with respect to the initial data and denoting the derivative by $\xi(s,t,y)$, we have
\begin{eqnarray}\label{3.20}
d \xi(s,t,y)&=&\lambda \nabla U(t,\gamma^{-1}(t,Y(s,t,y)))\nabla\gamma^{-1}(t,Y(s,t,y))\xi(s,t,y)dt\nonumber\\&&+\nabla^2
U(t,\gamma^{-1}(t,Y(s,t,y)))\nabla\gamma^{-1}(t,Y(s,t,y))\xi(s,t,y)d B(t)\nonumber\\
&=&:A_1(t,Y(s,t,y))\xi(s,t,y)dt+A_2(t,Y(s,t,y))\xi(s,t,y)d B(t),
\end{eqnarray}
with $\xi(s,t,y)|_{t=s}=I$.

Notice that the equation  (\ref{3.20}) is a linear  SDE with  bounded coefficients $A_{1}$\,, $A_{2}$ depending on the   process $Y(s,t,y)$,  by the Cauchy--Lipschitz theorem,  there is a  unique strong solution for equation (\ref{3.20})\,.
Moreover for  $d=1$, the unique strong solution is represented by
\begin{eqnarray}\label{3.21}
\xi(s,t,y)=\exp\Bigg(\int\limits_s^t[A_1(r,Y(s,r,y))-\frac{1}{2}A_2^2(r,Y(s,r,y))]dr +\int\limits_s^t
A_2(r,Y(s,r,y))d B(r)\Bigg).
\end{eqnarray}
By (\ref{3.10}) and Remark \ref{rem2.1}, $\nabla^2 U(t,x)$ is uniformly continuous in $x$,  then $A_1(s,r,Y(s,\cdot,\cdot))$ and $A_2(s,r,Y(s,\cdot,\cdot))$ are continuous in $[s,T]\times\mR^d$ almost surely.  Then for $d=1$\,,  the process
$\xi(s,t,y)$ is continuous in $(s,t,y)$ almost surely, and for every $p\geq 1$ there is a constant $C(p,T)>0$ such that
\begin{eqnarray}\label{3.22}
\sup_{y\in\mR^d}\sup_{0\leq s\leq
T}{\mathbb E}\sup_{s\leq t \leq T}\|\xi(s,t,y)\|^p\leq C(p,T).
\end{eqnarray}
For general $d>1$, the unique strong solution also has an obvious representation which is analogue of (\ref{3.21}), then using the same argument,   the continuity of $\xi(s,t,y)$ in $(s,t,y)$ and the inequality (\ref{3.22}) also hold. Then using the relationship between $X$ and $Y$, and (\ref{3.18}), (\ref{3.4}) holds true.

For every $x,y\in\mR^d$,  write $Y(s,t,x)$, $Y(s,t,y)$, $\xi(s,t,x)$, $\xi(s,t,y)$, $Y(s,t,x)-Y(s,t,y)$, $\xi(s,t,x)-\xi(s,t,y)$, $U(t,\gamma^{-1}(t,Y(s,t,x)))$ and $U(t,\gamma^{-1}(t,Y(s,t,y)))$ by $Y(x)$, $Y(y)$, $\xi(x)$, $\xi(y)$, $Y(x,y)$, $\xi(x,y)$, $U(\gamma^{-1}(Y(x)))$ and $U(\gamma^{-1}(Y(y)))$, respectively.  Then for every $q\geq 2$, by using It\^o's  formula
\begin{eqnarray}\label{3.23}
&&d \|\xi(x,y)\|^q\nonumber\\&=&q\lambda\|\xi(x,y)\|^{q-2}\langle \xi(x,y),
\nabla U(\gamma^{-1}(Y(x)))\nabla\gamma^{-1}(Y(x))\xi(x)-\nabla U(\gamma^{-1}(Y(y)))\nabla\gamma^{-1}(Y(y))\xi(y)\rangle dt\nonumber\\&&+
\frac{1}{2}q(q-1)\|\xi(x,y)\|^{q-2}tr([\nabla^2 U(\gamma^{-1}(Y(x)))\nabla\gamma^{-1}(Y(x))\xi(x)-\nabla^2 U(\gamma^{-1}(Y(y)))\nabla\gamma^{-1}(Y(y))
\nonumber\\&&\quad\times\xi(y)][\nabla^2 U(\gamma^{-1}(Y(x)))\nabla\gamma^{-1}(Y(x))\xi(x)-\nabla^2 U(\gamma^{-1}(Y(y)))\nabla\gamma^{-1}(Y(y))\xi(y)]^\top)dt
\nonumber\\&&+q\|\xi(x,y)\|^{q-2}\langle \xi(x,y),
[\nabla^2 U(\gamma^{-1}(Y(x)))\nabla\gamma^{-1}(Y(x))\xi(x)\nonumber\\&&\quad-\nabla^2 U(\gamma^{-1}(Y(y)))\nabla\gamma^{-1}(Y(y))\xi(y)]d B(t)\rangle
\nonumber\\&\leq& C(q)\Big[ \|\xi(x,y)\|^{q-1}\|\nabla U(\gamma^{-1}(Y(x)))\nabla\gamma^{-1}(Y(x))\xi(x)-\nabla U(\gamma^{-1}(Y(y)))\nabla\gamma^{-1}(Y(y))\xi(y)\|\nonumber\\&&+\|\xi(x,y)\|^{q-2}\|\nabla^2 U(\gamma^{-1}(Y(x)))\nabla\gamma^{-1}(Y(x))\xi(x)-\nabla^2 U(\gamma^{-1}(Y(y)))\nabla\gamma^{-1}(Y(y))\xi(y)\|^2\Big]dt
\nonumber\\&&+q\|\xi(x,y)\|^{q-2}\langle \xi(x,y),
[\nabla^2 U(\gamma^{-1}(Y(x)))\nabla\gamma^{-1}(Y(x))\xi(x)\nonumber\\&&\quad-\nabla^2 U(\gamma^{-1}(Y(y)))\nabla\gamma^{-1}(Y(y))\xi(y)]d B(t)\rangle.
\end{eqnarray}
Notice that
\begin{eqnarray}\label{3.24}
&&\|\nabla U(\gamma^{-1}(Y(x)))\nabla\gamma^{-1}(Y(x))\xi(x)-\nabla U(\gamma^{-1}(Y(y)))\nabla\gamma^{-1}(Y(y))\xi(y)\|\nonumber\\&\leq &
\|\nabla U(\gamma^{-1}(Y(x)))\nabla\gamma^{-1}(Y(x))\xi(x)-\nabla U(\gamma^{-1}(Y(y)))\nabla\gamma^{-1}(Y(x))\xi(x)\|\nonumber\\&&+\|\nabla U(\gamma^{-1}(Y(y)))\nabla\gamma^{-1}(Y(x))\xi(x)-\nabla U(\gamma^{-1}(Y(y)))\nabla\gamma^{-1}(Y(y))\xi(x)\|\nonumber\\&&+\|\nabla U(\gamma^{-1}(Y(y)))\nabla\gamma^{-1}(Y(y))\xi(x)
-\nabla U(\gamma^{-1}(Y(y)))\nabla\gamma^{-1}(Y(y))\xi(y)\|
\nonumber\\&\leq & \|\nabla^2 U\|_{L^\infty([0,T];\cC_b(\mR^d;\mR^{d\times d}\otimes\mR^d))}\|\nabla \gamma^{-1}\|^2_{L^\infty([0,T];\cC_b(\mR^d;\mR^{d\times d}))}|Y(x,y)|\|\xi(x)\|
\nonumber\\&&+ \|\nabla U\|_{L^\infty([0,T];\cC_b(\mR^d;\mR^{d\times d}))}\|\nabla^2 \gamma^{-1}\|_{L^\infty([0,T];\cC_b(\mR^d;\mR^{d\times d}\otimes\mR^d))}|Y(x,y)|\|\xi(x)\|
\nonumber\\&&+ \|\nabla U\|_{L^\infty([0,T];\cC_b(\mR^d;\mR^{d\times d}))}\|\nabla \gamma^{-1}\|_{L^\infty([0,T];\cC_b(\mR^d;\mR^{d\times d}))}\|\xi(x,y)\|
\nonumber\\&\leq & C\Big[|Y(x,y)|\|\xi(x)\|+\|\xi(x,y)\|\Big]
\end{eqnarray}
and
\begin{eqnarray}\label{3.25}
&&\|\nabla^2 U(\gamma^{-1}(Y(x)))\nabla\gamma^{-1}(Y(x))\xi(x)-\nabla^2 U(\gamma^{-1}(Y(y)))\nabla\gamma^{-1}(Y(y))\xi(y)\|
\nonumber\\&\leq&\|\nabla^2 U(\gamma^{-1}(Y(x)))\nabla\gamma^{-1}(Y(x))\xi(x)-\nabla^2 U(\gamma^{-1}(Y(y)))\nabla\gamma^{-1}(Y(x))\xi(x)\|\nonumber\\&&+\|\nabla^2 U(\gamma^{-1}(Y(y)))\nabla\gamma^{-1}(Y(x))\xi(x)-\nabla^2 U(\gamma^{-1}(Y(y)))\nabla\gamma^{-1}(Y(y))\xi(x)\|\nonumber\\&&+\|\nabla^2 U(\gamma^{-1}(Y(y)))\nabla\gamma^{-1}(Y(y))\xi(x)
-\nabla^2 U(\gamma^{-1}(Y(y)))\nabla\gamma^{-1}(Y(y))\xi(y)\|\nonumber\\&\leq & \|\nabla^2 U(\gamma^{-1}(Y(x)))-\nabla^2 U(\gamma^{-1}(Y(y)))\|\|\nabla \gamma^{-1}\|_{L^\infty([0,T];\cC_b(\mR^d;\mR^{d\times d}))}\|\xi(x)\|
\nonumber\\&&+ \|\nabla^2 U\|_{L^\infty([0,T];\cC_b(\mR^d;\mR^{d\times d}\otimes\mR^d))}\|\nabla^2 \gamma^{-1}\|_{L^\infty([0,T];\cC_b(\mR^d;\mR^{d\times d}\otimes\mR^d))}|Y(x,y)|\|\xi(x)\|
\nonumber\\&&+ \|\nabla^2 U\|_{L^\infty([0,T];\cC_b(\mR^d;\mR^{d\times d}\otimes\mR^d))}\|\nabla \gamma^{-1}\|_{L^\infty([0,T];\cC_b(\mR^d;\mR^{d\times d}))}\|\xi(x,y)\|
\nonumber\\&\leq & C\Big[\|\nabla^2 U(\gamma^{-1}(Y(x)))-\nabla^2 U(\gamma^{-1}(Y(y)))\| \|\xi(x)\| +|Y(x,y)|\|\xi(x)\|+\|\xi(x,y)\|\Big],
\end{eqnarray}
then for $t\in [s,T]$
\begin{eqnarray}\label{3.26}
&&\mE \|\xi(x,y)\|^q(t)\nonumber\\ &\leq& C(q,d)\mE\int\limits_s^t\|\xi(x,y)\|^q(r)dr+ C(q,d)\mE\int\limits_s^t\|\xi(x,y)\|^{q-1}(r)|Y(x,y)|(r)\|\xi(x)\|(r)dr\nonumber\\&& +C(q,d)\mE\int\limits_s^t\|\xi(x,y)\|^{q-2}(r)|Y(x,y)|^2(r)\|\xi(x)\|^2(r)dr \nonumber\\&&
+C(q,d)\mE\int\limits_s^t\|\xi(x,y)\|^{q-2}(r)\|\nabla^2 U(\gamma^{-1}(Y(x)))-\nabla^2 U(\gamma^{-1}(Y(y)))\|^2 \|\xi(x)\|^2dr
\nonumber\\&\leq&  C(q,d)\mE\int\limits_s^t\|\xi(x,y)\|^q(r)dr+ C(q,d)\mE\int\limits_s^t|Y(x,y)|^q(r)\|\xi(x)\|^q(r)dr\nonumber\\&&
+C(q,d)\mE\int\limits_s^t\|\nabla^2 U(\gamma^{-1}(Y(x)))-\nabla^2 U(\gamma^{-1}(Y(y)))\|^q \|\xi(x)\|^qdr
\nonumber\\&\leq&  C(q,d)\mE\int\limits_s^t\|\xi(x,y)\|^q(r)dr+ C(q,d)\Bigg[\mE\int\limits_s^t|Y(x,y)|^{2q}(r)dr\Bigg]^{\frac12}
\Bigg[\mE\int\limits_s^t|\xi(x)\|^{2q}(r)dr\Bigg]^{\frac12}\nonumber\\&&
+C(q,d)\Bigg[\mE\int\limits_s^t\|\nabla^2 U(\gamma^{-1}(Y(x)))-\nabla^2 U(\gamma^{-1}(Y(y)))\|^{2q} dr\Bigg]^{\frac12}\Bigg[\mE\int\limits_s^t\|\xi(x)\|^{2q}dr\Bigg]^{\frac12}.
\end{eqnarray}
Moreover  by (\ref{3.16}) and (\ref{3.22})
\begin{eqnarray}\label{3.27}
&&\mE \|\xi(x,y)\|^q(t)\nonumber\\ &\leq& C(q,d,T)\Bigg[|x-y|^q+\Bigg(\mE\int\limits_s^t\|\nabla^2 U(\gamma^{-1}(Y(x)))-\nabla^2 U(\gamma^{-1}(Y(y)))\|^{2q}dr\Bigg)^{\frac12}\Bigg].
\end{eqnarray}

Further  for every $p\geq 2$,  by the Burkholder--Davis--Gundy inequality
\begin{eqnarray}\label{3.28}
&&\sup_{0\leq s\leq
T}\mE \sup_{s\leq t\leq T}\|\xi(x,y)\|^p\nonumber\\ &\leq&
C(p,d,T)|x-y|^p \nonumber\\ && +
C(p,d,T)\mE\Bigg[\int\limits_0^T \|\nabla^2 U(\gamma^{-1}(Y(x)))-\nabla^2 U(\gamma^{-1}(Y(y)))\|^{4p}(r)dr\Bigg]^{\frac14}.
\end{eqnarray}

For given $p\geq2$, by the assumption in Theorem \ref{3.1}, there is a small enough real number $\delta=\delta(p)$ such that $F_\delta^{4p}(\cdot)$ (see (\ref{3.3})) is increasing and concave. On the other hand, since $U\in L^\infty([0,T];\cC^2_b({\mathbb R}^d;{\mathbb R}^d))$ and(\ref{3.10}) holds, for every $x,y\in\mR^d$, and the given real number $\delta(p)>0$, there is a constant $C(\delta)>0$ such that
\begin{eqnarray*}
\sup_{0\leq t\leq T}\|\nabla^2 U(t,x)-\nabla^2 U(t,y)\| \leq C(\delta)F_\delta(|x-y|) 1_{|x-y|<\delta} +C(\delta)|x-y|1_{|x-y|\geq \delta}\,.
\end{eqnarray*}
Then  by  (\ref{3.14}) and the fact that $F_\delta^{4p}(\cdot)$ is increasing
\begin{eqnarray}\label{3.29}
&&\sup_{0\leq s\leq
T}\mE \sup_{s\leq t\leq T}\|\xi(x,y)\|^p\nonumber\\ &\leq&
C|x-y|^p +
C\mE\Bigg[\int\limits_0^T F_\delta^{4p}(|\gamma^{-1}(Y(x))-\gamma^{-1}(Y(y))|) 1_{|\gamma^{-1}(Y(x))-\gamma^{-1}(Y(y))|<\delta}dt\Bigg]^{\frac14}
\nonumber\\ &&+C\mE\Bigg[\int\limits_0^T |\gamma^{-1}(Y(x))-\gamma^{-1}(Y(y))|^{4p} 1_{|\gamma^{-1}(Y(x))-\gamma^{-1}(Y(y))|\geq\delta}dt\Bigg]^{\frac14}
\nonumber\\ &\leq&
C|x-y|^p +C
\Bigg[\int\limits_0^T\mE F_\delta^{4p}(|\gamma^{-1}(Y(x))-\gamma^{-1}(Y(y))|1_{|\gamma^{-1}(Y(x))-\gamma^{-1}(Y(y))|<\delta}) dt\Bigg]^{\frac14}\nonumber\\&&+C\Bigg[\int\limits_0^T \mE|\gamma^{-1}(Y(x))-\gamma^{-1}(Y(y))|^{4p} dt\Bigg]^{\frac14}
\nonumber\\ &\leq&
C|x-y|^p +C
\Bigg[\int\limits_0^T F_\delta^{4p}(\mE[|\gamma^{-1}(Y(x))-\gamma^{-1}(Y(y))|1_{|\gamma^{-1}(Y(x))-\gamma^{-1}(Y(y))|<\delta}]) dt\Bigg]^{\frac14}\nonumber\\&\leq&C|x-y|^p +C F_\delta^{p}(2C_1|x-y|),
\end{eqnarray}
where $C_1=C(p,T)\vee 1$ and $C(p,T)$ is given in (\ref{3.17}).

Now by  $X(s,t,x)=\gamma^{-1}(t,Y(s,t,\gamma(s,x)))$, for every $p\geq 2$ and every $x,y\in\mR^d$ \begin{eqnarray*}
&&\|\nabla X(s,t,x)-\nabla X(s,t,y)\|\nonumber\\&=&\|\nabla \gamma^{-1}(t,Y(s,t,\gamma(s,x))) \nabla Y(s,t,\gamma(s,x))  \nabla \gamma(s,x)
\nonumber\\&&-\nabla \gamma^{-1}(t,Y(s,t,\gamma(s,y))) \nabla Y(s,t,\gamma(s,y))  \nabla \gamma(s,y)\|
\nonumber\\&=&\|[\nabla \gamma^{-1}(t,Y(s,t,\gamma(s,x))) -\nabla \gamma^{-1}(t,Y(s,t,\gamma(s,y))) ]\nabla Y(s,t,\gamma(s,x))  \nabla \gamma(s,x)
\nonumber\\&&+ \nabla \gamma^{-1}(t,Y(s,t,\gamma(s,y))) [\nabla Y(s,t,\gamma(s,x))-\nabla Y(s,t,\gamma(s,y)) ] \nabla \gamma(s,x)
\nonumber\\&&+  \nabla \gamma^{-1}(t,Y(s,t,\gamma(s,y))) \nabla Y(s,t,\gamma(s,y))[ \nabla \gamma(s,x)-\nabla \gamma(s,y)]\|
\nonumber\\&\leq&
\|\nabla\gamma^{-1}\|_{L^\infty([0,T];\cC_b^1(\mR^d;\mR^{d\times d}))}^2  \|Y(s,t,\gamma(s,x)) -Y(s,t,\gamma(s,y))\|\|\nabla Y(s,t,\gamma(s,x))\|
\nonumber\\&&+ \|\nabla \gamma^{-1}\|_{L^\infty([0,T];\cC_b(\mR^d;\mR^{d\times d}))}^2\|\nabla Y(s,t,\gamma(s,x))-\nabla Y(s,t,\gamma(s,y))\|
\nonumber\\&&+  \|\nabla \gamma^{-1}\|_{L^\infty([0,T];\cC_b^1(\mR^d;\mR^{d\times d}))} \|\nabla^2 \gamma\|_{L^\infty([0,T];\cC_b^1(\mR^d;\mR^{d\times d}\otimes\mR^d))}\| \nabla Y(s,t,\gamma(s,y))\| |x-y|
\nonumber\\&\leq&
C\Big[\|Y(s,t,\gamma(s,x)) -Y(s,t,\gamma(s,y))\|\|\nabla Y(s,t,\gamma(s,x))\|
\nonumber\\&&+ \|\nabla Y(s,t,\gamma(s,x))-\nabla Y(s,t,\gamma(s,y))\|+ \| \nabla Y(s,t,\gamma(s,y))\| |x-y|\Big].
\end{eqnarray*}
Thanks to the H\"{o}lder inequality, (\ref{3.22}) and (\ref{3.29}), for every $x,y\in\mR^d$, we have that
\begin{eqnarray}\label{3.30}
&& \sup_{0\leq s\leq
T}{\mathbb E}[\sup_{s\leq t \leq T}\|\nabla X(s,t,x)-\nabla X(s,t,y)\|^p]
\nonumber\\&\leq&
C\sup_{0\leq s\leq
T}\Big({\mathbb E}[\sup_{s\leq t \leq T}\|Y(s,t,\gamma(s,x)) -Y(s,t,\gamma(s,y))\|^{2p}]\Big)^{\frac12}
\Big({\mathbb E}\|\nabla Y(s,t,\gamma(s,x))\|^{2p}\Big)^{\frac12}
\nonumber\\&&+
 \sup_{0\leq s\leq
T}{\mathbb E}[\sup_{s\leq t \leq T}\|\nabla Y(s,t,\gamma(s,x))-\nabla Y(s,t,\gamma(s,y))\|^p]+C|x-y|^p
\nonumber\\&\leq&C|x-y|^p+CF_\delta^{p}(2C_1|\gamma(s,x)-\gamma(s,y)|)
\nonumber\\&\leq&C|x-y|^p+CF_\delta^{p}(3C_1|x-y|),
\end{eqnarray}
where in the last inequality we have used $\|\nabla \gamma\|_{L^\infty([0,T];{\mathcal C}_b({\mathbb R}^d;\mR^{d\times d}))}\leq 3/2$.
Then  for every $x$, $y\in\mR^d$ such that $0<|x-y|<\delta/(3C_1) \ (<r_0/(3C_1))$, by similar  calculations for (\ref{2.28})
\begin{eqnarray}\label{3.31}
&&\sup_{0\leq s\leq
T}{\mathbb E}[\sup_{s\leq t \leq T}\|\nabla X(s,t,x)-\nabla X(s,t,y)\|^p]
\nonumber \\
&\leq &C\Bigg[\int\limits_{r\leq 3C_1|x-y|}\frac{\phi(r)}{r}dr+ \phi(|x-y|)+|x-y|\int\limits_{3C_1|x-y|<r\leq \delta}\frac{\phi(r)}{r^2}dr|+|x-y|\Bigg]^p+C|x-y|^p
\nonumber \\
&\leq &C\Bigg[\int\limits_{r\leq |x-y|}\frac{\phi(r)}{r}dr+ \phi(|x-y|)+|x-y|\int\limits_{|x-y|<r\leq r_0}\frac{\phi(r)}{r^2}dr|\Bigg]^p+C|x-y|^p.
\end{eqnarray}
On account of (\ref{3.31}) and (\ref{3.4}), by using the H\"{o}lder inequality, we gain (\ref{3.6}).

(ii) Let $U^n$ be the unique solution of the parabolic problem (\ref{3.9}) associated to $b^n$, i.e.  $b$ is replaced by $b^n$  in equation (\ref{3.9}). By Theorem \ref{the2.1} (ii), $U^n\in L^\infty([0,T];\cC^2_b({\mathbb R}^d;\mR^{d}))\cap W^{1,\infty}([0,T];\cC_b({\mathbb R}^d;\mR^{d}))$ and
\begin{eqnarray}\label{3.32}
\lim_{n\rightarrow\infty}\| U^n-U\|_{L^\infty([0,T];\cC_b^2({\mathbb R}^d;\mR^{d}))}=0.
\end{eqnarray}
We set $\gamma_n(t,x):=x+U^n(t,x)$, then if $\lambda\geq \lambda_0$, $\{\gamma_n(t,x):=x+U^n(t,x)\}_n$ form nonsingular diffeomorphisms of class $\mathcal{C}^2$ uniformly in $t\in [0,T]$ and $n$. Moreover,
\begin{eqnarray}\label{3.33}
\lim_{n\rightarrow\infty}\| \gamma_n-\gamma\|_{L^\infty([0,T];\cC_b^2({\mathbb R}^d;\mR^d))}=\lim_{n\rightarrow\infty}\| \gamma_n^{-1}-\gamma^{-1}\|_{L^\infty([0,T];\cC_b^2({\mathbb R}^d;\mR^d))}=0.
\end{eqnarray}
To prove (\ref{3.7}) and (\ref{3.8}), it is sufficient to show  for every $p\geq 2$
\begin{eqnarray}\label{3.34}
\lim _{n\rightarrow\infty}\sup_{y\in {\mathbb R}^d} \sup_{0\leq s\leq
T}{\mathbb E}[\sup_{s\leq t \leq T}|Y^n(s,t,y)-Y(s,t,y)|^p]=0
\end{eqnarray}
and
\begin{eqnarray}\label{3.35}
\lim_{n\rightarrow \infty}\sup_{y\in {\mathbb R}^d} \sup_{0\leq s\leq
T}{\mathbb E}[\sup_{s\leq t \leq T}\|\nabla Y^n(s,t,y)-\nabla Y(s,t,y)\|^p]=0,
\end{eqnarray}
where $Y^n$ satisfies
\begin{eqnarray}\label{3.36}
\left\{
  \begin{array}{ll}
d Y^n(s,t)=\lambda U^n(t,\gamma^{-1}_n(t,Y^n(s,t)))dt+ [I+\nabla
U^n(t,\gamma^{-1}_n(t,Y^n(s,t)))] dB(t), \ t\in(s,T], \\
Y^n(s,t)|_{t=s}=y.
  \end{array}
\right.
\end{eqnarray}
For simplicity, we write $Y(s,t)$, $Y^n(s,t)$, $U(t,\gamma^{-1}(t,Y(s,t)))$ and $U^n(t,\gamma^{-1}_n(t,Y^n(s,t)))$ by $Y$, $Y^n$, $U(\gamma^{-1}(Y))$ and $U^n(\gamma^{-1}_n(Y^n))$, respectively. For every $q\geq 2$, by using the It\^{o} formula to $|Y^n-Y|^q$
\begin{eqnarray}\label{3.37}
&&d |Y^n-Y|^q\nonumber\\&=&q|Y^n-Y|^{q-2}\langle Y^n-Y,
U^n(\gamma^{-1}_n(Y^n))-U(\gamma^{-1}(Y))\rangle dt\nonumber\\&&+
\frac{1}{2}q(q-1)|Y^n-Y|^{q-2}tr([\nabla U^n(\gamma^{-1}_n(Y^n))-\nabla U(\gamma^{-1}(Y))]
\nonumber\\&&\quad\times[U^n(\gamma^{-1}_n(Y^n))-\nabla U(\gamma^{-1}(Y))]^\top)dt
\nonumber\\&&+q|Y^n-Y|^{q-2}\langle Y^n-Y,
[\nabla U^n(\gamma^{-1}_n(Y^n))-\nabla U(\gamma^{-1}(Y))]d B(t)\rangle
\nonumber\\&\leq& C(q)\Big[ |Y^n-Y|^{q-1}\|U^n(\gamma^{-1}_n(Y^n))-U(\gamma^{-1}(Y))\|\nonumber\\&&+|Y^n-Y|^{q-2}\|\nabla U^n(\gamma^{-1}_n(Y^n))-\nabla U(\gamma^{-1}(Y))\|^2\Big]dt
\nonumber\\&&+q|Y^n-Y|^{q-2}\langle Y^n-Y,
[\nabla U^n(\gamma^{-1}_n(Y^n))-\nabla U(\gamma^{-1}(Y))]d B(t)\rangle.
\end{eqnarray}
Then for every $t\in [s,T]$
\begin{eqnarray}\label{3.38}
&&\mE |Y^n-Y|^q(t)\nonumber\\ &\leq& C(q)\mE\int\limits_s^t |Y^n-Y|^{q-1}(r)\|U^n(\gamma^{-1}_n(Y^n))-U(\gamma^{-1}(Y))\|(r)dr \nonumber\\&& +C(q)\mE\int\limits_s^t|Y^n-Y|^{q-2}(r)\|\nabla U^n(\gamma^{-1}_n(Y^n))-\nabla U(\gamma^{-1}(Y))\|^2(r)dr.
\end{eqnarray}
First we have
\begin{eqnarray}\label{3.39}
&&\|U^n(\gamma^{-1}_n(Y^n))-U(\gamma^{-1}(Y))\|\nonumber\\&=&\|U^n(\gamma^{-1}_n(Y^n))-U(\gamma^{-1}_n(Y^n))
+U(\gamma^{-1}_n(Y^n))-U(\gamma^{-1}(Y^n))+U(\gamma^{-1}(Y^n))
-U(\gamma^{-1}(Y))\|
\nonumber\\&\leq & \|U^n-U\|_{L^\infty([0,T];\cC_b(\mR^d;\mR^{d}))}
+\|\nabla U\|_{L^\infty([0,T];\cC_b(\mR^d;\mR^{d\times d}))}\|\gamma^{-1}_n-\gamma^{-1}\|_{L^\infty([0,T];\cC_b(\mR^d;\mR^{d}))}\nonumber\\&& +\|\nabla U\|_{L^\infty([0,T];\cC_b(\mR^d;\mR^{d\times d}))}\|\nabla \gamma^{-1}\|_{L^\infty([0,T];\cC_b(\mR^d;\mR^{d\times d}))}|Y^n-Y|
\nonumber\\&\leq & C\Big[\|U^n-U\|_{L^\infty([0,T];\cC_b(\mR^d;\mR^{d}))}+
\|\gamma^{-1}_n-\gamma^{-1}\|_{L^\infty([0,T];\cC_b(\mR^d;\mR^{d}))}
+|Y^n-Y|\Big]
\end{eqnarray}
and similarly
\begin{eqnarray}\label{3.40}
&&\|\nabla U^n(\gamma^{-1}_n(Y^n))-\nabla U(\gamma^{-1}(Y))\|^2\nonumber\\&\leq & C\Big[\|\nabla U^n-\nabla U\|_{L^\infty([0,T];\cC_b(\mR^d;\mR^{d\times d}))}\nonumber\\&&+\|\nabla^2 U\|_{L^\infty([0,T];\cC_b(\mR^d;\mR^{d\times d}\otimes\mR^d))}\|\gamma^{-1}_n-\gamma^{-1}\|_{L^\infty([0,T];\cC_b(\mR^d;\mR^{d}))}
\nonumber\\&&+\|\nabla^2 U\|_{L^\infty([0,T];\cC_b(\mR^d;\mR^{d\times d}\otimes\mR^d))}\|\nabla \gamma^{-1}\|_{L^\infty([0,T];\cC_b(\mR^d;\mR^{d\times d}))}|Y^n-Y|\Big]^2
\nonumber\\&\leq & C\Big[\|\nabla U^n-\nabla U\|^2_{L^\infty([0,T];\cC_b(\mR^d;\mR^{d\times d}))}+\|\gamma^{-1}_n-\gamma^{-1}\|^2_{L^\infty([0,T];\cC_b(\mR^d;\mR^{d}))}
+|Y^n-Y|^2\Big].
\end{eqnarray}
Since $Y_n$ and $Y$ satisfy (\ref{3.36}) and (\ref{3.20}), respectively, and the coefficients are globally Lipschitz continuous, for every $q\geq 2$, there is a positive constant $C(q)$ such that
\begin{eqnarray}\label{3.41}
\sup_n\sup_{y\in {\mathbb R}^d} \sup_{0\leq s\leq
T}\sup_{s\leq t \leq T}{\mathbb E}[|Y^n(s,t,y)-Y(s,t,y)|^q]\leq C(q).
\end{eqnarray}
Then by (\ref{3.38})--(\ref{3.41}) and a Gr\"onwall type argument, we obtain
\begin{eqnarray}\label{3.42}
&&\sup_{y\in {\mathbb R}^d} \sup_{0\leq s\leq
T}\sup_{s\leq t\leq T}\mE |Y^n(s,t,y)-Y(s,t,y)|^q\nonumber\\ &\leq& C(q,T)\Big[\|U^n-U\|^q_{L^\infty([0,T];\cC_b^1(\mR^d;\mR^{d}))}+
\|\gamma^{-1}_n-\gamma^{-1}\|^q_{L^\infty([0,T];\cC_b(\mR^d;\mR^{d}))}
\Big].
\end{eqnarray}

Now for every $p\geq 2$,  by (\ref{3.37}) and the Burkholder--Davis--Gundy inequality
\begin{eqnarray}\label{3.43}
&&\mE \sup_{s\leq t \leq T}|Y^n(s,t,y)-Y(s,t,y)|^p\nonumber\\ &\leq& C(p)\mE\int\limits_s^T |Y^n-Y|^{p-1}(r)\|U^n(\gamma^{-1}_n(Y^n))-U(\gamma^{-1}(Y))\|(r)dr \nonumber\\&& +C(p)\mE\int\limits_s^T|Y^n-Y|^{p-2}(r)\|\nabla U^n(\gamma^{-1}_n(Y^n))-\nabla U(\gamma^{-1}(Y))\|^2(r)dr
 \nonumber\\&&+C(p)\mE\Bigg[\int\limits_s^T|Y^n-Y|^{2p-2}(r)\|\nabla U^n(\gamma^{-1}_n(Y^n))-\nabla U(\gamma^{-1}(Y))\|^2(r)dr\Bigg]^{\frac12}.
\end{eqnarray}
Then thanks to (\ref{3.39}), (\ref{3.40}), (\ref{3.42}) and the Gr\"onwall inequality,
\begin{eqnarray}\label{3.44}
&&\sup_{y\in {\mathbb R}^d} \sup_{0\leq s\leq
T}\mE \sup_{s\leq t\leq T}|Y^n(s,t,y)-Y(s,t,y)|^p\nonumber\\ &\leq& C(p,T)\Bigg[\|U^n-U\|^{p}_{L^\infty([0,T];\cC_b^1(\mR^d;\mR^{d}))}+
\|\gamma^{-1}_n-\gamma^{-1}\|^{p}_{L^\infty([0,T];\cC_b(\mR^d;\mR^{d}))}
\Bigg],
\end{eqnarray}
which implies (\ref{3.34}) by  (\ref{3.32})  and (\ref{3.33})\,.

Differentiate $Y^n(s,t,y)$  with respect to initial data~(denoted by $\xi^n(s,t,y)$ or $\xi^n$), then~$\xi^{n}$ satisfies equation~(\ref{3.20}) with $b$ replaced by $b^{n}$\,,  and
by the boundedness of the coefficients and initial value~($I$)\,,
 for every $q\geq 2$, there is a positive real number $C(q)$ such that
\begin{eqnarray}\label{3.45}
\sup_n\sup_{y\in {\mathbb R}^d} \sup_{0\leq s\leq
T}\sup_{s\leq t \leq T}{\mathbb E}[\|\xi^n(s,t,y)\|^q]\leq C(q).
\end{eqnarray}
%
%
%
Moreover observing that
\begin{eqnarray}\label{3.46}
&&\|\nabla U^n(\gamma^{-1}_n(Y^n))\nabla\gamma^{-1}_n(Y^n)\xi^n-\nabla U(\gamma^{-1}(Y))\nabla\gamma^{-1}(Y)\xi\|\nonumber\\&=&
\|\nabla U^n(\gamma^{-1}_n(Y^n))\nabla\gamma^{-1}_n(Y^n)\xi^n-\nabla U(\gamma^{-1}_n(Y^n))\nabla\gamma^{-1}_n(Y^n)\xi^n
\nonumber\\
&&+\nabla U(\gamma^{-1}_n(Y^n))\nabla\gamma^{-1}_n(Y^n)\xi^n-\nabla U(\gamma^{-1}(Y^n))\nabla\gamma^{-1}_n(Y^n)\xi^n
\nonumber\\
&&+
\nabla U(\gamma^{-1}(Y^n))\nabla\gamma^{-1}_n(Y^n)\xi^n-\nabla U(\gamma^{-1}(Y))\nabla\gamma^{-1}_n(Y^n)\xi^n
\nonumber\\
&&+\nabla U(\gamma^{-1}(Y))\nabla\gamma^{-1}_n(Y^n)\xi^n-\nabla U(\gamma^{-1}(Y))\nabla\gamma^{-1}(Y^n)\xi^n
\nonumber\\
&&+\nabla U(\gamma^{-1}(Y))\nabla\gamma^{-1}(Y^n)\xi^n-\nabla U(\gamma^{-1}(Y))\nabla\gamma^{-1}(Y)\xi^n
\nonumber\\
&&+
\nabla U(\gamma^{-1}(Y))\nabla\gamma^{-1}(Y)\xi^n-\nabla U(\gamma^{-1}(Y))\nabla\gamma^{-1}(Y)\xi\|
\nonumber\\&\leq & \|\nabla U^n-\nabla
U\|_{L^\infty([0,T];\cC_b(\mR^d;\mR^{d\times d}))}\|\nabla\gamma^{-1}_n\|_{L^\infty([0,T];\cC_b(\mR^d;\mR^{d\times d}))}\|\xi^n\|
\nonumber\\
&&+\|\nabla^2 U\|_{L^\infty([0,T];\cC_b(\mR^d;\mR^{d\times d}\otimes\mR^d))}\|\gamma^{-1}_n-\gamma^{-1}\|_{L^\infty([0,T];\cC_b(\mR^d;\mR^d))}
\|\nabla\gamma^{-1}_n\|_{L^\infty([0,T];\cC_b(\mR^d;\mR^{d\times d}))}\|\xi^n\|
\nonumber\\&&
+\|\nabla^2 U\|_{L^\infty([0,T];\cC_b(\mR^d;\mR^{d\times d}\otimes\mR^d))}\|\nabla \gamma^{-1}\|_{L^\infty([0,T];\cC_b(\mR^d;\mR^{d\times d}))}
\nonumber\\ &&\quad
\times |Y^n-Y|
\|\nabla\gamma^{-1}_n\|_{L^\infty([0,T];\cC_b(\mR^d;\mR^{d\times d}))}\|\xi^n\|
\nonumber\\
&&+\|\nabla U\|_{L^\infty([0,T];\cC_b(\mR^d;\mR^{d\times d}))}\|\nabla\gamma^{-1}_n-\nabla\gamma^{-1}\|_{L^\infty([0,T];\cC_b(\mR^d;\mR^{d\times d}))}\|\xi^n\|
\nonumber\\
&&+\|\nabla U\|_{L^\infty([0,T];\cC_b(\mR^d;\mR^{d\times d}))}\|\nabla^2\gamma^{-1}\|_{L^\infty([0,T];\cC_b(\mR^d;\mR^{d\times d}\otimes\mR^d))}|Y^n-Y|\|\xi^n\|
\nonumber\\
&&+\|\nabla U\|_{L^\infty([0,T];\cC_b(\mR^d;\mR^{d\times d}))}\|\nabla\gamma^{-1}\|_{L^\infty([0,T];\cC_b(\mR^d;\mR^{d\times d}))}\|\xi^n-\xi\|
\nonumber\\&\leq &
C\Big[
 \|\nabla U^n-\nabla
U\|_{L^\infty([0,T];\cC_b(\mR^d;\mR^{d\times d}))}
+\|\gamma^{-1}_n-\gamma^{-1}\|_{L^\infty([0,T];\cC_b^1(\mR^d;\mR^{d}))}+|Y^n-Y|\Big]\|\xi^n\|
\nonumber\\&&
+C\|\xi^n-\xi\|
\end{eqnarray}
and by the boundedness of $\|\nabla^2 U \|$
\begin{eqnarray}\label{3.47}
&&\|\nabla^2 U^n(\gamma^{-1}_n(Y^n))\nabla\gamma^{-1}_n(Y^n)\xi^n-\nabla^2
U(\gamma^{-1}(Y))\nabla\gamma^{-1}(Y)\xi\|
\nonumber\\&=&\|\nabla^2 U^n(\gamma^{-1}_n(Y^n))\nabla\gamma^{-1}_n(Y^n)\xi^n-
\nabla^2 U(\gamma^{-1}_n(Y^n))\nabla\gamma^{-1}_n(Y^n)\xi^n
\nonumber\\&&+ \nabla^2 U(\gamma^{-1}_n(Y^n))\nabla\gamma^{-1}_n(Y^n)\xi^n-\nabla^2 U(\gamma^{-1}(Y))\nabla\gamma^{-1}_n(Y^n)\xi^n
\nonumber\\&&+
\nabla^2 U(\gamma^{-1}(Y))\nabla\gamma^{-1}_n(Y^n)\xi^n-\nabla^2 U(\gamma^{-1}(Y))\nabla\gamma^{-1}(Y^n)\xi^n
\nonumber\\&&+
\nabla^2 U(\gamma^{-1}(Y))\nabla\gamma^{-1}(Y^n)\xi^n-\nabla^2 U(\gamma^{-1}(Y))\nabla\gamma^{-1}(Y)\xi^n
\nonumber\\&&+
\nabla^2 U(\gamma^{-1}(Y))\nabla\gamma^{-1}(Y)\xi^n-\nabla^2
U(\gamma^{-1}(Y))\nabla\gamma^{-1}(Y)\xi\|
\nonumber\\&\leq & C\Big[\|\nabla^2 U^n-\nabla^2 U\|_{L^\infty([0,T];\cC_b(\mR^d;\mR^{d\times d}\otimes\mR^d))} +\|\nabla \gamma^{-1}_n-\nabla\gamma^{-1}\|_{L^\infty([0,T];\cC_b(\mR^d;\mR^{d\times d}))}
\nonumber\\&&+|Y^n-Y| \Big]\|\xi^n\|+C\|\nabla^2 U(\gamma^{-1}_n(Y^n))-\nabla^2 U(\gamma^{-1}(Y))\|\|\xi^n\|
+C\|\xi^n-\xi\|\,,
\end{eqnarray}
 then by similar calculations for  (\ref{3.37})--(\ref{3.44})\,, we have
 \begin{eqnarray}\label{3.48}
&&\sup_{y\in {\mathbb R}^d} \sup_{0\leq s\leq
T}\mE \sup_{s\leq t\leq T}\|\xi^n(s,t,y)-\xi(s,t,y)\|^p\nonumber\\ &\leq&
C(p,T)\Big[\|U^n-U\|^{p}_{L^\infty([0,T];\cC_b^2(\mR^d;\mR^{d}))}+
\|\gamma^{-1}_n-\gamma^{-1}\|^{p}_{L^\infty([0,T];\cC_b^1(\mR^d;\mR^{d}))}
\Big] \nonumber\\&& +
C(p,T)\sup_{y\in {\mathbb R}^d} \sup_{0\leq s\leq
T}\mE\Big[\int\limits_s^T \|\nabla^2 U(\gamma^{-1}_n(Y^n))-\nabla^2 U(\gamma^{-1}(Y))\|^{2p}\|\xi^n\|^{2p}dr\Big]^{\frac12}
\nonumber\\ &\leq&
C(p,T)\Big[\|U^n-U\|^{p}_{L^\infty([0,T];\cC_b^2(\mR^d;\mR^{d}))}+
\|\gamma^{-1}_n-\gamma^{-1}\|^{p}_{L^\infty([0,T];\cC_b^1(\mR^d;\mR^{d}))}
\Big] \nonumber\\&& +
C(p,T)\sup_{y\in {\mathbb R}^d}\mE\Big[\int\limits_0^T \|\nabla^2 U(\gamma^{-1}_n(Y^n))-\nabla^2 U(\gamma^{-1}(Y))\|^{4p}dr\Big]^{\frac14}\,.
\end{eqnarray}

Further we apply same calculations for (\ref{3.29}) to $\|\nabla^2 U(\gamma^{-1}_n(Y^n))-\nabla^2 U(\gamma^{-1}(Y))\|$ to get
\begin{eqnarray}\label{3.49}
&&\sup_{y\in {\mathbb R}^d}\mE\Bigg[\int\limits_0^T \|\nabla^2 U(\gamma^{-1}_n(Y^n(t,y)))-\nabla^2 U(\gamma^{-1}(Y(t,y)))\|^{4p}dr\Bigg]^{\frac14}\nonumber\\ &\leq&
C(p,T,\delta)\sup_{y\in {\mathbb R}^d}\mE\Bigg[\int\limits_0^T F_\delta^{4p}(|\gamma^{-1}_n(Y(t,y))-\gamma^{-1}(Y(t,y))|) 1_{|\gamma^{-1}_n(Y(t,y))-\gamma^{-1}(Y(t,y))|<\delta}dt\Bigg]^{\frac14}
\nonumber\\ &&+C(p,T,\delta)\sup_{y\in {\mathbb R}^d}\mE\Bigg[\int\limits_0^T |\gamma^{-1}_n(Y(t,y))-\gamma^{-1}(Y(t,y))|^{4p} 1_{|\gamma^{-1}_n(Y(t,y))-\gamma^{-1}(Y(t,y))|\geq\delta}dt\Bigg]^{\frac14}
\nonumber\\ &\leq&
C\Bigg[\int\limits_0^T F_\delta^{4p}(\sup_{y\in {\mathbb R}^d}\mE|\gamma^{-1}_n(Y(t,y))-\gamma^{-1}(Y(t,y))|) dt\Bigg]^{\frac14}+C
\|\gamma^{-1}_n-\gamma^{-1}\|^{p}_{L^\infty([0,T];\cC_b(\mR^d;\mR^{d}))}
\nonumber\\ &\leq&
C\|\gamma^{-1}_n-\gamma^{-1}\|^{p}_{L^\infty([0,T];\cC_b(\mR^d;\mR^{d}))}
\nonumber\\&&+C F_\delta^{p}\Big(C\sup_{y\in {\mathbb R}^d} \sup_{0\leq s\leq
T}\mE\sup_{s\leq r\leq
T}|Y^n(s,r,y)-Y(s,r,y)|+\|\gamma^{-1}_n-\gamma^{-1}\|_{L^\infty([0,T];\cC_b(\mR^d;\mR^{d}))}\Big)
\nonumber\\ &\leq&
C\|\gamma^{-1}_n-\gamma^{-1}\|^{p}_{L^\infty([0,T];\cC_b(\mR^d;\mR^{d}))}
\nonumber\\&&+C F_\delta^{p}\Big(C\|U^n-U\|_{L^\infty([0,T];\cC_b^1(\mR^d;\mR^{d}))}+C\|\gamma^{-1}_n-
\gamma^{-1}\|_{L^\infty([0,T];\cC_b(\mR^d;\mR^{d}))}\Big),
\end{eqnarray}
where in the last inequality we  use (\ref{3.44}), and in the third inequality we use
\begin{eqnarray*}
|\gamma_n^{-1}(Y^n)-\gamma^{-1}(Y)|&\leq& \|\nabla\gamma^{-1}_n\|_{L^\infty([0,T];\cC_b(\mR^d;\mR^{d\times d}))}|Y^n-Y|+\|\gamma^{-1}_n-
\gamma^{-1}\|_{L^\infty([0,T];\cC_b(\mR^d;\mR^d))}\\ &\leq&C|Y^n-Y|+\|\gamma^{-1}_n-
\gamma^{-1}\|_{L^\infty([0,T];\cC_b(\mR^d;\mR^d))}.
\end{eqnarray*}
Now from  (\ref{3.48}) and (\ref{3.49}),  (\ref{3.35}) holds by (\ref{3.32}) and (\ref{3.33})\,.   $\Box$

\begin{remark} \label{rem3.1} Let $b^n$ be given in Theorem \ref{the3.1} (ii). By Liouville's theorem,  since $b^n$ is smooth,  we have the following Euler identity
\begin{eqnarray*}
\det(\nabla_xX^n(s,t,x))=\exp\Big(\int\limits^t_s\div  b^n(r,X^n(s,r,x))dr\Big), \ \ 0\leq s\leq t \leq T,
\end{eqnarray*}
where $\det(\cdot)$ denotes the determinant of a matrix. In view of (\ref{3.7}) and (\ref{3.8}), if  $b\in L^\infty([0,T]; {\mathcal C}_b({\mathbb R}^d ;{\mathbb R}^d))$ such that (\ref{3.2}) holds, and $\div  b\in L^1([0,T];L^1_{loc}({\mathbb R}^d))$\,, up to choosing a subsequence, one derives
\begin{eqnarray}\label{3.50}
\det(\nabla_xX(s,t,x))=\exp\Big(\int\limits^t_s\div  b(r,X(s,r,x))dr\Big), \quad 0\leq s\leq t \leq T.
\end{eqnarray}
\end{remark}

Since the inverse of $X$ ($X^n$): $X^{-1}$ ($X^{-1}_n$) satisfies an equation which has the same form as the original one except the drift has opposite sign, by Theorem \ref{the3.1} we have
\begin{corollary} \label{cor3.1} Let $b,b^n$, $X(s,t,x)$ and $X^n(s,t,x)$ be stated in Theorem \ref{the3.1}. For every $p\geq 1$, there is a constant $C(p,d,T)>0$ such that
\begin{eqnarray}\label{3.51}
\sup_{0\leq s\leq
T}{\mathbb E}\sup_{s\leq t \leq T}|X^{-1}(s,t,x)|^p+\sup_{x\in\mR^d}\sup_{0\leq s\leq
T}{\mathbb E}\sup_{s\leq t \leq T}\|\nabla X^{-1}(s,t,x)\|^p\leq C(p,d,T)
\end{eqnarray}
and
\begin{eqnarray}\label{3.52}
&&\lim _{n\rightarrow\infty}\sup_{x\in {\mathbb R}^d} \sup_{0\leq s\leq
T}{\mathbb E}\sup_{s\leq t \leq T}[|X^{-1}_n(s,t,x)-X^{-1}(s,t,x)|]^p
\nonumber\\ &=&\lim _{n\rightarrow\infty}\sup_{x\in {\mathbb R}^d} \sup_{0\leq s\leq
T}{\mathbb E}\sup_{s\leq t \leq T}[\|\nabla X^{-1}_n(s,t,x)-\nabla X^{-1}(s,t,x)\|]^p=0.
\end{eqnarray}
Moreover, for every $x,y\in \mR^d$, there is another constant $C(p,d,T)>0$ such that
\begin{eqnarray}\label{3.53}
\sup_{0\leq s\leq
T}{\mathbb E}[\sup_{s\leq t \leq T}|X^{-1}(s,t,x)-X^{-1}(s,t,y)|^p]\leq C(p,d,T)|x-y|^p
\end{eqnarray}
and
\begin{eqnarray}\label{3.54}
&&\sup_{0\leq s\leq
T}{\mathbb E}[\sup_{s\leq t \leq T}\|\nabla X^{-1}(s,t,x)-\nabla X^{-1}(s,t,y)\|^p]
\nonumber \\
&\leq &C\Bigg[\int\limits_{r\leq |x-y|}\frac{\phi(r)}{r}dr+ \phi(|x-y|)+|x-y|\int\limits_{|x-y|<r\leq r_0}\frac{\phi(r)}{r^2}dr|\Bigg]^p1_{|x-y|<r_0}+C|x-y|^p.
\end{eqnarray}
\end{corollary}

\begin{remark} \label{rem3.2}
(i) Let $r_0=1/2$ and $\phi(r)=C|\log(r)|^{-\alpha}$ for $r\in (0,r_0)$ with some $\alpha>1$. From Remark \ref{rem2.3} (ii)
\begin{eqnarray*}
\int\limits_{r\leq |x-y|}\frac{\phi(r)}{r}dr+ \phi(|x-y|)+|x-y|\int\limits_{|x-y|<r\leq r_0}\frac{\phi(r)}{r^2}dr+|x-y|\leq\frac{C(\alpha)}{|\log(|x-y|)|^{\alpha-1}}.
\end{eqnarray*}
Let $F_\delta$ be defined by (\ref{3.3}). We get
\begin{eqnarray}\label{3.55}
F_\delta(r)\leq \frac{C}{|\log(r)|^{\alpha-1}}=:\tilde{F}_\delta(r), \quad  r\in [0,\delta], \quad \delta<\frac12.
\end{eqnarray}
Notice that  given $p\geq 1$,  if $r \in (0,\delta)$ and $\delta<\exp(p-p\alpha-1)$,
\begin{eqnarray}\label{3.56}
\frac{d}{dr} \frac{1}{|\log(r)|^{p(\alpha-1)}}=\frac{p(\alpha-1)}{r|\log(r)|^{p(\alpha-1)+1}}\geq 0
\end{eqnarray}
and
\begin{eqnarray}\label{3.57}
\frac{d^2}{dr^2} \frac{1}{|\log(r)|^{p(\alpha-1)}}=\frac{p(\alpha-1)}{|\log(r)|^{p(\alpha-1)+1}r^2}
 \Bigg[\frac{p\alpha-p+1}{|\log(r)|}-1\Bigg]\leq 0\,,
\end{eqnarray}
then the function $\tilde{F}_\delta^p(r)=C^p|\log(r)|^{-p(\alpha-1)}$ is increasing and concave (we define $|\log(0)|^{-1}=0$) on~$[0,\delta]$. Thus the assumptions in Theorem \ref{3.1} hold. Theorem \ref{3.1} is applicable mutatis mutandis.
\end{remark}

From Remark \ref{rem3.2}, we draw the following result.
\begin{corollary} \label{cor3.2} Let $b\in L^\infty([0,T];\cC_b(\mR^d;\mR^d))$. Suppose that there are two real numbers $C>0$  and $\alpha>1$ such that for every $x\in \mR^d$
\begin{eqnarray}\label{3.58}
|b(t,x)-b(t,y)|\leq \frac{C}{|\log(|x-y|)|^\alpha}, \quad  {\rm for\; all}  \ \ y\in B_{\frac12}(x), \ t\in [0,T].
\end{eqnarray}

(i) SDE
(\ref{3.1}) has a unique continuous adapted solution $\{X(s,t,x), \
t\in [s,T ], \ \omega \in\Omega\}$, which forms a stochastic diffeomorphisms flow.

(ii) Estimates (\ref{3.4})--(\ref{3.5}) hold for $X$ and denote its inverse by $X^{-1}$, then (\ref{3.51}) and (\ref{3.53}) hold as well. Moreover, for every $x,y\in\mR^d$
\begin{eqnarray}\label{3.59}
&&\sup_{0\leq s\leq
T}{\mathbb E}\sup_{s\leq t \leq T}[\|\nabla X(s,t,x)-\nabla X(s,t,y)\|+\|\nabla X^{-1}(s,t,x)-\nabla X^{-1}(s,t,y)\|]^p
\nonumber \\ &\leq&
C(p,T)\Bigg[\frac{1_{|x-y|< \frac12}}{|\log(|x-y|)|^{p(\alpha-1)}}+|x-y|^p1_{|x-y|\geq \frac12}\Bigg].
\end{eqnarray}
\end{corollary}

\section{Stochastic transport equations}\label{sec4}
\setcounter{equation}{0}
\begin{definition} \label{def4.1}
Let $b\in L^1([0,T];L^1_{loc}({\mathbb R}^d;{\mathbb R}^d))$ such that $\div  b\in L^1([0,T];L^1_{loc}({\mathbb R}^d))$, and let  $u_0\in L^\infty({\mathbb R}^d)$.
A stochastic field $u$ is called a weak $L^\infty$-solution of (\ref{1.1}) if $u\in L^\infty(\Omega\times[0,T];L^\infty({\mathbb R}^d))$ and for every
$\varphi\in{\mathcal C}_0^\infty({\mathbb R}^d)$, $\int_{{\mathbb R}^d}\varphi(x)u(t,x)dx$
has a continuous modification which is an ${\mathcal F}_t$-semimartingale and
for every  $t\in [0,T]$
\begin{eqnarray}\label{4.1}
\int\limits_{{\mathbb R}^d}\varphi(x)u(t,x)dx&=&\int\limits_{{\mathbb R}^d}\varphi(x)u_0(x)dx+
\int\limits^t_0\int\limits_{{\mathbb R}^d}\div (b(s,x)\varphi(x))u(s,x)dxds\nonumber\\&&
+\sum_{i=1}^d\int\limits^t_0\circ dB_i(s)\int\limits_{{\mathbb R}^d}\partial_{x_i}\varphi(x)u(s,x)dx,  \quad
{\mathbb P}-a.s..
\end{eqnarray}
\end{definition}
Then we  state our main result.
\begin{theorem} \label{the4.1} \textbf{(Existence and uniqueness)} Let $d\geq 1$. Suppose $b\in L^\infty([0,T];\cC_b(\mR^d;\mR^d))$ such that (\ref{3.2}) and (\ref{3.3}) hold. Further suppose  that $\div b\in L^1([0,T];L^1_{loc}(\mR))$ for $d=1$ or there exists $q\in (2,+\infty)$ such that
\begin{eqnarray}\label{4.2}
\div b\in L^q([0,T]\times {\mathbb R}^d), \quad  d\geq 1.
\end{eqnarray}
Then there exists a unique weak $L^\infty$-solution to the Cauchy problem (\ref{1.1}).  Moreover, the unique
weak solution can be represented by $u(t,x)=u_0(X^{-1}(t,x))$, with $X(t,x)$ being the unique
strong solution of the associated stochastic differential equation  (\ref{3.1}) with $s=0$.
 \end{theorem}

\noindent
\textbf{Proof.} First, we prove that $u(t,x)=u_0(X^{-1}(t,x))$ is a weak $L^\infty$-solution of (\ref{1.1}). Let $b^n$ and $X^n$ be given in Theorem \ref{the3.1}, and let $X^{-1}_n$ be the inverse of $X^n$. Since $b^n$ is smooth,  $u^n=u_0(X^{-1}_n)$ is the unique weak $L^\infty$-solution of the following Cauchy problem~(\cite[Theorems 16 and 20]{FGP1})
\begin{eqnarray}\label{4.3}
\left\{
  \begin{array}{ll}
\partial_tu^n(t,x)+b^n(t,x)\cdot\nabla u^n(t,x)
+\sum_{i=1}^d\partial_{x_i}u^n(t,x)\circ\dot{B}_i(t)=0, \ (t,x)\in(0,T)\times {\mathbb R}^d, \\
u^n(t,x)|_{t=0}=u_0(x), \  x\in{\mathbb R}^d.
  \end{array}
\right.
\end{eqnarray}
Then for every
$\varphi\in{\mathcal C}_0^\infty({\mathbb R}^d)$ and
 every  $t\in [0,T]$
\begin{eqnarray}\label{4.4}
\int\limits_{{\mathbb R}^d}\varphi(x)u_0(X^{-1}_n(t,x))dx&=&\int\limits_{{\mathbb R}^d}\varphi(x)u_0(x)dx+
\int\limits^t_0\int\limits_{{\mathbb R}^d}\div (b^n(s,x)\varphi(x))u_0(X^{-1}_n(s,x))dxds\nonumber\\&&
+\sum_{i=1}^d\int\limits^t_0\circ dB_i(s)\int\limits_{{\mathbb R}^d}\partial_{x_i}\varphi(x)u_0(X^{-1}_n(s,x))dx
\nonumber\\&=&\int\limits_{{\mathbb R}^d}\varphi(x)u_0(x)dx+
\int\limits^t_0\int\limits_{{\mathbb R}^d}\div (b^n(s,x)\varphi(x))u_0(X^{-1}_n(s,x))dxds\nonumber\\&&
+\sum_{i=1}^d\int\limits^t_0dB_i(s)\int\limits_{{\mathbb R}^d}\partial_{x_i}\varphi(x)u_0(X^{-1}_n(s,x))dx
\nonumber\\ && +\frac12 \int\limits^t_0 ds\int\limits_{{\mathbb R}^d}\Delta\varphi(x)u_0(X^{-1}_n(s,x))dx  \quad
{\mathbb P}-a.s.,
\end{eqnarray}
where in the last inequality we have used the relationship between the Stratonovich integral and the It\^{o} integral.

By (\ref{2.32}) and (\ref{4.2})
\begin{eqnarray}\label{4.5}
\lim_{n\rightarrow \infty}\|b^n-b\|_{L^\infty([0,T];\cC_b(\mR^d;\mR^d))}=0
\end{eqnarray}
and
\begin{eqnarray}\label{4.6}
\lim_{n\rightarrow \infty}\|\div(b^n\varphi)-\div(b\varphi)\|_{L^1([0,T];L^1(\mR^d))}=0.
\end{eqnarray}
Now  thanks to Corollary \ref{cor3.1}, by taking $n$ to infinity in (\ref{4.4}), (\ref{4.1}) holds for $u(t,x)=u_0(X^{-1}(t,x))$.

It remains to check the uniqueness and observing that the equation is linear, it suffices to prove that $u\equiv 0$ a.s. if the initial data vanishes. We only check the uniqueness for $d>1$, the case $d=1$ being similar and easier. Let $\varrho_n$ be given by (\ref{2.8}) and  set $u_n=u\ast \varrho_n$, then
\begin{eqnarray}\label{4.7}
\partial_tu_n(t,x)+b(t,x)\cdot\nabla u_n(t,x)
+\sum_{i=1}^d\partial_{x_i}u_n(t,x)\circ\dot{B}_i(t)=e_n(t,x),
\end{eqnarray}
with
\begin{eqnarray}\label{4.8}
e_n(t,x)=b(t,x)\cdot\nabla u_n(t,x)-(b\cdot\nabla u)\ast \varrho_n(t,x).
\end{eqnarray}
For for every $t>0$\,,  the  It\^{o}'s formula yields that
\begin{eqnarray*}
u_n(t,X(t,x))=\int_0^te_n(s,X(s,x))ds,
\end{eqnarray*}
which implies for every $\varphi\in \mathcal{C}_0^\infty({\mathbb R}^d)$ and almost all $\omega\in\Omega$
\begin{eqnarray}\label{4.9}
&&\int\limits_{{\mathbb R}^d}u_n(t,X(t,x))\varphi(x)dx
\nonumber\\&=&
\int\limits^t_0\int\limits_{{\mathbb R}^d}e_n(s,X(s,x))\varphi(x)dx
\nonumber\\&=&\int\limits^t_0\int\limits_{{\mathbb R}^d}e_n(s,x)\varphi(X^{-1}(s,x))
\det(\nabla_xX^{-1}(s,x))dxds
\nonumber\\&=&-\int\limits^t_0\int\limits_{{\mathbb R}^d}\int\limits_{{\mathbb R}^d}[\div b(s,x)-\div b(s,y)]\varphi(X^{-1}(s,x))
\det(\nabla_xX^{-1}(s,x))u(s,y)\varrho_n(x-y)dydxds
\nonumber\\&&+\int\limits^t_0\int\limits_{{\mathbb R}^d}\int\limits_{{\mathbb R}^d}[b(s,x)-b(s,y)]\cdot \nabla_x(\varphi(X^{-1}(s,x))
\det(\nabla_xX^{-1}(s,x)))u(s,y)\varrho_n(x-y)dydxds
\nonumber\\&=&: I_1^n(t)+I_2^n(t).
\end{eqnarray}

In view of Theorem \ref{the3.1}, $\varphi(X^{-1}(s,x))\det(\nabla_xX^{-1}(s,x))$ is continuous in $(s,x)$. We claim that for almost all $\omega\in\Omega$, $\varphi(X^{-1}(s,\cdot))$  has a compact support in $x$ uniformly in $s$. In fact,   without loss of generality we assume there is a real number $R>0$ such that the
support of $\varphi$ is in $B_R$, if the assertion is false, there is a sequence $\{(t_k,x_k)\}_{k\geq1}\subset [0,t]\times B_R$ such that $\lim_{k\rightarrow \infty}|X(\omega,t_k,x_k)|=+\infty$. But on the other hand, by Theorem \ref{the3.1}, for every $k\geq 1$
\begin{eqnarray*}
|X(\omega,t_k,x_k)|&\leq& |X(\omega,t_k,x_k)-X(\omega,t_k,x_1)|+|X(\omega,t_k,x_1)|\\ &\leq& C|x_k-x_1|+ |X(\omega,t_k,x_1)|\\ &\leq& 2CR+\sup_{0\leq s\leq T}|X(\omega,s,x_1)\leq C.
\end{eqnarray*}
Therefore, by (\ref{4.2}), taking  $n$ to infinity yields   $I_1^n(t)\rightarrow 0$\,,  $\mP-$a.s..

Noticing that
\begin{eqnarray*}
&&\nabla_x(\varphi(X^{-1}(s,x)))
\det(\nabla_xX^{-1}(s,x))\nonumber\\&=&\nabla_{X^{-1}}\varphi(X^{-1}(s,x))\nabla_xX^{-1}(s,x)
\det(\nabla_xX^{-1}(s,x))+\varphi(X^{-1}(s,x))\nabla_x(
\det(\nabla_xX^{-1}(s,x))),
\end{eqnarray*}
if for almost all $\omega\in\Omega$, $\nabla_x(
\det(\nabla_xX^{-1}(\cdot,\cdot)))\in L^1([0,T];L^1_{loc}(\mR^d))$, by  the dominated convergence theorem,  $I_2^n(t)\rightarrow 0$, $\mP$-a.s.. Now let us check $\nabla_x(
\det(\nabla_xX^{-1}(\cdot,\cdot)))\in L^1([0,T];L^1_{loc}(\mR^d))$ and it is equivalent to show that $\nabla_x(
\det(\nabla_xX(\cdot,\cdot)))\in L^1([0,T];L^1_{loc}(\mR^d))$.
From (\ref{3.50})
\begin{eqnarray*}
\nabla_x\det(\nabla_xX(t,x))=\det(\nabla_xX(t,x))\nabla_x\int\limits^t_0\div  b(s,X(s,x))ds.
\end{eqnarray*}
By Theorem \ref{the3.1}, $\det(\nabla_xX(s,x))$ is continuous in $(s,x)$, we need to show \begin{eqnarray}\label{4.10}
\nabla_x\int\limits^t_0\div  b(s,X(s,x))ds\in L^1([0,T];L^1_{loc}(\mR^d)), \quad \mP-a.s..
\end{eqnarray}

To do this we consider the following backward parabolic equation
\begin{eqnarray}\label{4.11}
\left\{
\begin{array}{ll}
\partial_{t}V(t,x) +\frac{1}{2}\Delta V(t,x)+b(t,x)\cdot \nabla V(t,x)=\div b(t,x), \quad (t,x)\in (0,T)\times {\mathbb R}^d, \\
V(T,x)=0, \quad  x\in{\mathbb R}^d.
  \end{array}
\right.
\end{eqnarray}
Noticing  $\div b\in L^q([0,T]\times {\mathbb R}^d)$, there is a unique $V\in L^q([0,T];W^{2,q}({\mathbb R}^d))\cap W^{1,q}([0,T];L^q({\mathbb R}^d))$ solving (\ref{4.11}), and there is a constant $C(q,T)>0$ such that
\begin{eqnarray}\label{4.12}
\|V\|_{L^q([0,T];W^{2,q}({\mathbb R}^d))}+\|\partial_tV\|_{L^q([0,T]\times {\mathbb R}^d)} \leq C(q,T)\|\div b\|_{L^q([0,T]\times {\mathbb R}^d)}.
\end{eqnarray}
Moreover, by Sobolev's imbedding theorem,
\begin{eqnarray}\label{4.13}
\sup_{0\leq t\leq T}\|V\|_{W^{1,q}({\mathbb R}^d)} \leq C(q,T)\|\div b\|_{L^q([0,T]\times {\mathbb R}^d)}.
\end{eqnarray}
Next we assume that for every $t\in [0,T]$, $V(t)$ is smooth, otherwise one can follow an approximation argument, then  the It\^o's formula yields
\begin{eqnarray*}
V(t,X(t,x))-V(0,x)-\int\limits_0^t\nabla V(s,X(s,x))dB(s)=\int\limits_0^t\div b(s,X(s,x))ds,
\end{eqnarray*}
which implies that
\begin{eqnarray}\label{4.14}
\nabla_x\int\limits^t_0\div  b(s,X(s,x))ds&=&\nabla_X V(t,X(t,x))\nabla_xX(t,x)-\nabla_xV(0,x)
\nonumber\\&&-\int\limits_0^t\nabla^2_X V(s,X(s,x))\nabla_xX(s,x)dB(s).
\end{eqnarray}
By Theorem \ref{the3.1}, (\ref{4.12}) and (\ref{4.13}), for almost all $\omega\in\Omega$, the first two terms in the righthand hand side in (\ref{4.14}) are in $L^1([0,T];L^1_{loc}(\mR^d))$. For every $R>0$
\begin{eqnarray}\label{4.15}
&&\mE\int\limits_0^T\int\limits_{|x|\leq R}\Bigg|\int\limits_0^t\nabla^2_X V(s,X(s,x))\nabla_xX(s,x)dB(s)\Bigg|^2dxdt\nonumber\\&=&\int\limits_0^T\int\limits_{|x|\leq R}\mE\int\limits_0^t|\nabla^2_X V(s,X(s,x))\nabla_xX(s,x)|^2dsdxdt
\nonumber\\&\leq &C\int\limits_{|x|\leq R}\int\limits_0^T\Big[\mE\|\nabla^2_X V(s,X(s,x))\|^q\Big]^{\frac{2}{q}}\Big[\mE\|\nabla_xX(s,x)\|^\frac{2q}{q-2}\Big]^{\frac{q-2}{q}}dsdx
\nonumber\\&\leq &C\Big[\sup_{x\in\mR^d}\mE\sup_{0\leq s \leq T}|\nabla_xX(s,x)|^\frac{2q}{q-2}\Big]^{\frac{q-2}{q}}
\Bigg[\int\limits_{|x|\leq R}\int\limits_0^T\Big[\mE|\nabla^2_X V(s,X(s,x))|^q\Big]dsdx\Bigg]^{\frac{2}{q}}
\nonumber\\&\leq &C
\Bigg[\mE\int\limits_0^T\int\limits_{|x|\leq R}|\nabla^2_X V(s,X(s,x))|^qdxds\Bigg]^{\frac{2}{q}}
\nonumber\\&\leq &C\Bigg[\sup_{x\in\mR^d}\mE\sup_{0\leq s \leq T}\|\nabla_xX^{-1}(s,x)\|^q\Bigg]^{\frac{2}{q}}
\Bigg[\int\limits_0^T\int\limits_{\mR^d}|\nabla^2_x V(s,x)|^qdxds\Bigg]^{\frac{2}{q}}<+\infty,
\end{eqnarray}
where in the fifth line we have used (\ref{3.4}), and in the last inequality we have used (\ref{3.51}) and (\ref{4.12}).

Then from (\ref{4.9})--(\ref{4.15}), by taking $n$ to infinity, we have
$\int_{{\mathbb R}^d}u(t,X(t,x))\varphi(x)dx=0$\,, that is $u(t,X(t,x))=0$ for almost everywhere $x\in \mR^d$ and almost all $\omega\in\Omega$. Because $X(t,x)$ is a stochastic diffeomorphisms flow associated with (\ref{3.1}) with $s=0$, we have $u(t,x)=0$ for almost everywhere $x\in \mR^d$ and almost all $\omega\in\Omega$.  The proof is complete. $\Box$

\begin{remark} \label{rem4.1} Without the stochastic perturbation, even if the drift is bounded and H\"{o}lder continuous, the deterministic equation
possesses multiple $L^\infty$-solutions (\cite[Section 6.1]{FGP1}). So the noise has a regularization effect.
\end{remark}

\begin{remark} \label{rem4.2}
(i)  SDE (\ref{3.1}) with  $s=0$ has a unique continuous adapted solution $\{X(t,x), \ t\in [0,T ], \ \omega \in\Omega\}$, which forms a stochastic flow of diffeomorphisms, and by Theorem \ref{the4.1} (\ref{1.1}) has a unique weak $L^\infty$-solution which can be wrote by $u_0(X^{-1}(t,x))$. Thus, if $u_0\in \cC_b(\mR^d)$, for almost all $\omega\in \Omega$, and every $t\in [0,T]$, $u(t,\cdot)$ is bounded and continuous.  Moreover, if $u_0\in W^{1,p}(\mR^d)$ with $p\in [1,+\infty]$, we have the following chain rule
\begin{eqnarray*}
\nabla_x(u_0(X^{-1}(t,x)))=\nabla_xu_0(X^{-1}(t,x))
\nabla_xX^{-1}(t,x),
\end{eqnarray*}
so $u(t,\cdot)\in W^{1,p}_{loc}(\mR^d)$ almost surely, for every $t\in [0,T]$.
However, for the deterministic equation, even if the
uniqueness is established, the persistence of the above properties (continuity and Sobolev differentiability) for solutions are missing~\cite{CLR2}.

(ii) We can further  establish the existence and uniqueness of $W^{1,p}$-solutions with $p\in [1,+\infty]$ as well if $b\in L^\infty([0,T];\cC_b(\mR^d;\mR^d))$ such that (\ref{3.2}) holds. More precisely there is a unique $u\in L^\infty(\Omega\times[0,T];L^p({\mathbb R}^d))$ such that

(1) for every
$\varphi\in{\mathcal C}_0^\infty({\mathbb R}^d)$, $\int_{{\mathbb R}^d}\varphi(x)u(t,x)dx$
has a continuous modification which is an ${\mathcal F}_t$-semimartingale and
for every  $t\in [0,T]$, (\ref{4.1}) holds;

(2) for almost all $\omega\in\Omega$, $\nabla u\in L^\infty([0,T];L^p_{loc}(\mR^d;\mR^d))$.

The proof is similar to that  of Theorem \ref{the4.1}~(\cite[Theorem 1.1]{WDGL}, \cite[Theorem 25]{FGP1}). However, if $d\geq 2$, the deterministic equation does not exist such strong solutions~(\cite[Theorem 1.2]{WDGL})\,, which means that  noise prevents the singularity for solutions.

(iii) There are also some results on $\cap _{p\geq 1}W^{1,p}_{loc}$ solutions~(\cite{FF2}) and $L^p$ solutions~(\cite{CO})\,.

\end{remark}

\section{Conclusions}\label{sec4}\setcounter{equation}{0}
Recently, there have been a broad research on the uniqueness of $L^\infty$ solutions
for the stochastic transport equation
\begin{eqnarray}\label{5.1}
\left\{
  \begin{array}{ll}
\partial_tu(t,x)+b(t,x)\cdot\nabla u(t,x)
+\sum_{i=1}^d\partial_{x_i}u(t,x)\circ\dot{B}_i(t)=0, \quad (t,x)\in(0,T)\times {\mathbb R}^d, \\
u(t,x)|_{t=0}=u_0(x), \quad  x\in{\mathbb R}^d,
  \end{array}
\right.
\end{eqnarray}
with non-Lipschitz drift. Most of these works are concentrated on the drift which is H\"{o}lder continuous in spatial variable uniformly in time. The question for the uniqueness  when $b$ is only bounded is still open. In this study, we established the existence and uniqueness of $L^\infty$ solutions only assuming $b$ is bounded and Dini continuous in spatial variable uniformly in time. Compared with the existing research, the result is new.

We solve the Cauchy problem (\ref{5.1}) by the method of stochastic characteristics. Therefore, we should prove the existence of stochastic diffeomorphisms flow for the following stochastic differential equation
\begin{eqnarray}\label{5.2}
dX(t)=b(t,X(t))dt+dB(t), \quad  t\in(0,T], \quad  X(t)|_{t=0}=x.
\end{eqnarray}
To reach the goal, we use the It\^{o}-Tanaka trick  to
transform the SDE (\ref{5.2}) with bounded and Dini continuous drift to an equivalent new SDE with Lipschitz coefficients via a non-singular diffeomorphism $\Phi(t,x)=x+U(t,x)$, where $U(T-t,x)=:V(t,x)$ satisfies a vector-valued parabolic partial differential equation of second order which has the form
\begin{eqnarray}\label{5.3}
\left\{\begin{array}{ll}
\partial_{t}V(t,x)=\frac{1}{2}\Delta  V(t,x)+b(t,x)\cdot \nabla V(t,x)
+b(t,x)-\lambda V(t,x), \ (t,x)\in (0,T)\times {\mathbb R}^d, \\
V(0,x)=0, \  x\in{\mathbb R}^d.  \end{array}\right.
\end{eqnarray}
There are two things we need to do. The first one is choosing a proper function space on which the It\^{o} formula is applicable, and the second one is the boundedness estimate for the gradient of $V$. We accomplish these issues by fetching $L^\infty(0,T;\cC^2_b(\mR^d))\cap W^{1,2}(0,T;\cC_b(\mR^d))$ as the workspace. When $b\in L^\infty(0,T;\mathcal{C}^{\alpha}_b(\mathbb{R}^d;\mathbb{R}^d))$,
these estimates for solutions have been established by
Flandoli, Gubinelli and Priola~\cite{FGP1}.
Noticing that, here we only assume that
$b$ is Dini continuous in $x$, so we should extend the Schauder theory for (\ref{5.3}) to $W^{2,\infty}$ theory. We accomplish these estimates in Section 2, and then establish the stochastic diffeomorphisms flow for (\ref{5.2}) in Section 3. These results are new as well.

We remark that the method used to establish $W^{2,\infty}$ estimates for (\ref{5.3}) can be applied to found the $W^{1,\infty}$ estimates for solutions of second order parabolic equations driven by Browian motion or general L\'{e}vy noise.
The $W^{1,p}$ ($p\in [2,\infty)$) theory, stochastic BMO estimates and Schauder theory have been founded by many researchers, see  \cite{Kry96,Kim15,DuL,HWW,WDL}. There are few works to deal the $W^{1,\infty}$ estimates. Therefore, the study of the $W^{1,\infty}$ property of solutions
to stochastic parabolic equations is of very high importance. For simplicity, and without loss of generality, here we only give a brief calculation for the $W^{1,\infty}$ estimate for heat equation driven by Brownian noise:
\begin{eqnarray}\label{5.4}
\left\{\begin{array}{ll}
du(t,x)=\frac{1}{2}\Delta u(t,x)dt+f(t,x)dB(t), \ (t,x)\in (0,T)\times {\mathbb R}^d, \\
u(0,x)=0, \  x\in{\mathbb R}^d, \end{array}\right.
\end{eqnarray}
where $f$ is bounded and Dini continuous in $x$ uniformly in $t$. From (\ref{5.4}), then
\begin{eqnarray}\label{5.5}
u(t,x)=\int\limits_0^tK(t-s,\cdot)\ast f(s,\cdot)(x)dB(s),  \quad \ (t,x)\in [0,T]\times\mR^d,
\end{eqnarray}
where $K(t,x)=(2\pi t)^{-\frac{d}{2}}e^{-\frac{|x|^2}{2t}}, \ t>0, \ x\in\mathbb{R}^d$. For $1\leq i\leq d$, we first differentiate $u$ in $x_i$, and then use the It\^{o} isometry, to get
\begin{eqnarray*}
&&\mE|\partial_{x_i}u(t,x)|^2\nonumber \\&=&\int\limits_0^t\Big|\int\limits_{\mR^d}\partial_{x_i}K(t-s,x-y) f(s,y)dy\Big|^2ds\nonumber \\ &=& \int\limits_0^t\Big|\int\limits_{|x-y|>(t-s)^\theta}\partial_{x_i}K(t-s,x-y) [f(s,y)-f(s,x)]dy\nonumber \\ && +\int\limits_{|x-y|\leq (t-s)^\theta}\partial_{x_i}K(t-s,x-y) [f(s,y)-f(s,x)]dy\Big|^2ds \nonumber \\ &\leq& 2\int\limits_0^t\Big|\int\limits_{|x-y|>(t-s)^\theta}\partial_{x_i}K(t-s,x-y) [f(s,y)-f(s,x)]dy\Big|^2ds\nonumber \\ && +2\int\limits_0^t\Big|\int\limits_{|x-y|\leq (t-s)^\theta}\partial_{x_i}K(t-s,x-y) [f(s,y)-f(s,x)]dy\Big|^2ds,
\end{eqnarray*}
where $\theta\in (0,1/2)$. Similar calculations from (\ref{2.12}) to (\ref{2.17}) used here again, we conclude that
\begin{eqnarray*}
\mE|\partial_{x_i}u(t,x)|^2<\infty.
\end{eqnarray*}
Moreover, we also get an analogue of (\ref{2.7}) for $\partial_{x_i}u$, i.e. for every $x,y\in \mR$ and every $t\in [0,T]$
\begin{eqnarray*}
&&\Big[\mE|\partial_{x_i}u(t,x)-\partial_{y_i}u(t,y)|^2\Big]^{\frac12} \nonumber \\ &\leq&
C\left[\int\limits_{r\leq |x-y|}\frac{\phi(r)}{r}dr+ \phi(|x-y|)+|x-y|\int\limits_{|x-y|<r\leq r_0}\frac{\phi(r)}{r^2}dr\right]1_{|x-y|<r_0}+C|x-y|.
\end{eqnarray*}
Moreover, if $f$ is H\"{o}lder-Dini or strong H\"{o}lder or weak H\"{o}lder continuous with the H\"{o}lder-Dini or strong H\"{o}lder or weak H\"{o}lder function $\phi$, then
\begin{eqnarray*}
\Big[\mE|\partial_{x_i}u(t,x)-\partial_{y_i}u(t,y)|^2\Big]^{\frac12} \leq C\phi(|x-y|).
\end{eqnarray*}

\section*{Acknowledgements}
This research was partly supported by the NSF of China grants 11501577, 11771123, 11771207.

\end{document}